\documentclass[a4paper,leqno,12pt]{amsart}
\usepackage[leqno]{amsmath}
\usepackage{amstext,amssymb,amsthm, enumitem}
\usepackage{pgfplots}
\usepackage{graphicx}
\usepackage{tikz}
\usepackage[cp1250]{inputenc}
\usepackage{float}
\usetikzlibrary{arrows}
\newtheorem{lemma}{Lemma}
\newtheorem{example}{Example}
\newtheorem*{lemma6}{Lemma 6$'$}
\newtheorem*{lemma7}{Lemma 7$'$}
\newtheorem*{theorem*}{Theorem}
\newtheorem{theorem}{Theorem}
\newtheorem{proposition}{Proposition}
\theoremstyle{definition}
\newtheorem*{prof*}{Proof}
\newtheorem*{proof*}{}

\numberwithin{equation}{section}

\tikzset{dots/.append style={ultra thick, fill=none}}

\begin{document}

\title[mapping properties of modified maximal operators]{On relations between weak and strong type inequalities for modified maximal operators on non-doubling metric measure spaces}

\author{Dariusz Kosz}
\address{ 
	\newline Faculty of Pure and Applied Mathematics
	\newline Wroc\l aw University of Science and Technology 
	\newline Wyb. Wyspia\'nskiego 27 
	\newline 50-370 Wroc\l aw, Poland
	\newline \textit{Dariusz.Kosz@pwr.edu.pl}	
}
\begin{abstract} In this article we investigate a special class of non-doubling metric measure spaces in order to describe the possible configurations of $P_{k,\rm s}^{\rm c}$, $P_{k,\rm s}$, $P_{k,\rm w}^{\rm c}$ and $P_{k,\rm w}$, the sets of all $p \in [1, \infty]$ for which the weak and strong type $(p,p)$ inequalities hold for the centered and non-centered modified Hardy--Littlewood maximal operators, $M^{\rm c}_k$ and $M_k$, $k \geq 1$. For any fixed $k$ we describe the necessary conditions that $P_{k,\rm s}^{\rm c}$, $P_{k,\rm s}$, $P_{k,\rm w}^{\rm c}$ and $P_{k,\rm w}$ must satisfy in general and illustrate each admissible configuration with a properly chosen non-doubling metric measure space. We also give some partial results related to an analogous problem stated for varying $k$.
	
	
\end{abstract}

\thanks{
	$\noindent \textit{Key words:}$ modified maximal operators; weak and strong type inequalities; non-doubling metric measure spaces.	
} 
\maketitle

\section{Introduction}

Let $\mathfrak{X} = (X, \rho, \mu)$ be a metric measure space with a metric $\rho$ and a Borel measure $\mu$. Throughout this article, unless otherwise stated, we assume that $(X, \rho)$ is bounded (that is, $\rm{diam}$$(X) = \sup\{\rho(x,y) \colon x, y \in X\} < \infty$) and the measure of each ball is finite and strictly positive. We also emphasize that each space we deal with later on is separable. By $B(x,r) = B_\rho(x,r)$ we denote the open ball centered at $x \in X$ with radius $r>0$. If we do not specify the center point and radius we write simply $B$. According to this notation, for a parameter $k \geq 1$, we define the $\textit{modified Hardy--Littlewood}$ $\textit{maximal operators}$, centered $M_k^{\rm c}$ and non-centered $M_k$, by
\begin{displaymath}
M_k^{\rm c}f(x) = \sup_{r > 0} \frac{1}{\mu(B(x,kr))} \int_{B(x,r)} |f| d\mu, \qquad x \in X,
\end{displaymath}
and
\begin{displaymath}
M_kf(x) = \sup_{B \ni x} \frac{1}{\mu(kB)} \int_B |f| d \mu , \qquad x \in X,
\end{displaymath}
respectively. Here $kB$ refers to the ball concentric with $B$ and
of radius $k$ times that of $B$. Note that, in general, neither the center nor the radius of a ball as a set are uniquely determined. Moreover, in the case $k > 1$, it is possible that for some $x, y \in X$ and $r,s > 0$ we have $B(x, r) = B(y, s)$, while $B(x, kr) \neq B(y, ks)$. If $k=1$, then the modified operators coincide with the standard Hardy--Littlewood maximal operators, non-centered and centered, and hence we will write shortly $M^{\rm c}$ or $M$ instead of $M_1^{\rm c}$ or $M_1$. Finally, let us make it clear that in the case of arbitrary $\mathfrak{X}$ the balls $B$ such that $|B| = 0$ or $|kB| = \infty$ are omitted in the definitions of $M_k^{\rm c}$ and $M_k$ (in the extreme case we use the convention that the supremum of the empty set is $0$). 

In this paper we investigate mapping properties of $M_k^{\rm c}$ and $M_k$ in the context of $L^p$ and weak $L^p$ function spaces for $p \in [1, \infty]$. So far, most of the work in this area was devoted to the case $p=1$, especially to study the weak type $(1,1)$ boundedness. There were several articles focused on the general description of the situations in which the weak type $(1,1)$ inequality must occur (see \cite{NTV}, \cite{Sa} and \cite{Te}, for example). Finally, it was proven in \cite{St1} that $M_k^{\rm c}$ and $M_k$ are of weak type $(1,1)$ for $k \geq 2$ and $k \geq 3$, respectively, in the case of any metric measure space with a measure that is finite on bounded sets. Moreover, it is known that the above ranges of the parameter $k$ are sharp in the sense that for any $k < 2$ (or $k < 3$) one can find a metric measure space such that $M_k^{\rm c}$ (or $M_k$) is not of weak type $(1,1)$. The suitable examples are given in \cite{Sa} and \cite{St2} (see also \cite{SS}, where certain details justifying the correctness of the construction described in \cite{Sa} are given). The aim of this article is to show as many as possible different admissible configurations of the sets of $p$ for which the weak and strong type $(p,p)$ inequalities hold, by using similar structures as those occuring in \cite{St2}. We study two cases, $k$ fixed or varying. 

Let us introduce the notation $A_1 \lesssim A_2$ (equivalently, $A_2 \gtrsim A_1$), which means that $A_1 \leq CA_2$ with a positive constant $C$ independent of significant quantities (in particular, $A_1 = \infty$ implies $A_2 = \infty$). We write $A_1 \simeq A_2$ if $A_1 \lesssim A_2$ and $A_2 \lesssim A_1$ hold simultaneously. Moreover, for a given measurable function $f \geq 0$ and a set $E \subset X$ of strictly positive measure we denote the average value of $f$ on $E$ by 
\begin{displaymath}
	A_E(f) = \frac{1}{\mu(E)}\int_{E} f \, d\mu. 
\end{displaymath}

Recall that for $p \in [1, \infty)$ the space $L^{p, \infty} = L^{p, \infty}(\mathfrak{X})$ consists of all measurable functions $g$ such that $\|g\|_{p, \infty} := \sup_{\lambda > 0}\{ \lambda \, \mu(E_\lambda(g))^{1/p} \} < \infty$, where $E_\lambda(g) := \{x \in X \colon |g(x)| > \lambda\}$ is the level set of $g$. Similarly, if $\|g\|_{p} := \big( \int_X |g|^p d\mu \big)^{1/p} < \infty$, then $g \in L^p = L^{p}(\mathfrak{X})$. Furthermore, the space $L^\infty = L^\infty(\mathfrak{X})$ is defined analogously by using $\|g\|_\infty := \inf\{C \geq 0 \colon |g| \leq C \textrm{ almost everywhere}\}$. Accordingly, we say that an operator $T$ is of strong (or weak) type $(p,p)$ for some $p \in [1, \infty]$, if $T$ is bounded on $L^p$ (or $T$ is bounded from $L^p$ to $L^{p,\infty}$), which means that $ \| T g \|_p \lesssim \| g \|_p$ (or $ \| T g \|_{p, \infty} \lesssim \| g \|_p$) holds uniformly in $f \in L^p$. Here we use the convention $L^{\infty, \infty}=L^\infty$. 

As a starting point of our considerations we explain a specific technique of combining different metric measure spaces, which will be often used later on. Fix $k_0 \geq 1$. Let $\Lambda$ be a (finite or not) set of positive integers and for each $n \in \Lambda$ consider a metric measure space $\mathfrak{X}_n = (X_n, \rho_n, \mu_n)$. We introduce $\rho_n'$ and $\mu_n'$ by rescaling (if necessary) $\rho_n$ and $\mu_n$, respectively, in such a way that $\rm{diam}\it(X_n)$ with respect to $\rho_n'$ does not exceed $1$ and $\mu_n'(X_n) \leq 2^{-n}$. Then, assuming that $X_{n_1} \cap X_{n_2} = \emptyset$ for any $n_1 \neq n_2$, we construct the space $\mathfrak{X} =  (X, \rho, \mu)$ as follows.
Denote $X = \bigcup_{n \in \Lambda} X_n$. We define the metric $\rho$ on $X$ by
\begin{displaymath}
\rho(x,y) = \left\{ \begin{array}{rl}
\rho_n'(x,y) & \textrm{if }  \{x,y\} \subset X_n \textrm{ for some } n \in \Lambda,   \\
k_0+1 & \textrm{otherwise,} \end{array} \right. 
\end{displaymath} 
and the measure $\mu$ on $X$ by
\begin{displaymath}
\mu(E) = \sum_{n \in \Lambda} \mu_n'(E \cap X_n), \qquad E \subset X.
\end{displaymath}
In the following proposition we describe some relations between the mapping properties of the maximal operators associated with $\mathfrak{X}$ and $\mathfrak{X}_n$, $n \in \Lambda$. 
\begin{proposition}
	Define $\mathfrak{X}$ as above for a fixed $k_0 \geq 1$. If $1 \leq k \leq k_0$, then the modified maximal operator $M_{k,\mathfrak{X}}$ (or $M_{k,\mathfrak{X}}^{\rm c}$) is of weak (respectively strong) type $(p,p)$ for some $p \in [1, \infty]$ if and only if for each $n \in \Lambda$ the operator $M_{k,\mathfrak{X}_n}$ (or $M_{k,\mathfrak{X}_n}^{\rm c}$) satisfies the weak (respectively strong) type $(p,p)$ inequality with a constant $\tilde{c} = \tilde{c}(k,p)$ that does not depend on $n$.   
\end{proposition}

\begin{prof*}
First note that the process of rescaling metrics and measures, which was used in the construction of $\mathfrak{X}$, does not affect the studied mapping properties of the associated maximal operators $M_{k,\mathfrak{X}_n}^{\rm c}$ and $M_{k,\mathfrak{X}_n}$, $n \in \Lambda$. Thus, without any loss of generality, we can simply assume that the spaces $\mathfrak{X}_n$ are the rescaled ones, that is, $\rm{diam}$$(X_n) \leq 1$ (with respect to $\rho_n$) and $\mu_n(X_n) \leq 2^{-n}$. Fix $1 \leq k \leq k_0$ and $p \geq 1$. To make the proof clear, assume that we study only the strong type $(p,p)$ of the non-centered operator (the other options can be considered similarly). Observe that if we take $f \in L^p(\mathfrak{X}_n)$ for some $n \in \Lambda$ and next we extend $f$ to $F \in L^p(\mathfrak{X})$, setting $F(x)=0$ for $ x \in X \setminus X_n$, then $\|F\|_p = \|f\|_p$ (here the symbol $\| \cdot \|_p$ refers to function spaces over different measure spaces) and $M_{k, \mathfrak{X}}(F)(x) = M_{k,\mathfrak{X}_n}(f)(x)$ for any $x \in X_n$. Hence, $\|M_{k, \mathfrak{X}_n}(f)\|_p / \|f\|_p \leq \|M_{k, \mathfrak{X}}(F)\|_p / \|F\|_p$. Thus, we conclude that if $\|M_{k, \mathfrak{X}}(g)\|_p \leq  \tilde{c}(k,p) \|g\|_p$ holds for every $g \in L^p(\mathfrak{X})$, then the operators $M_{k,\mathfrak{X}_n}$, $n \in \Lambda$, satisfy the adequate inequalities with the same constant $\tilde{c}(k,p)$. Now assume that each operator $M_{k,\mathfrak{X}_n}$, $n \in \Lambda$, satisfies the strong type $(p,p)$ inequality with a constant $\tilde{c} = \tilde{c}(k,p)$. Let $f \in L^p(\mathfrak{X})$ and define $f_n \in L^p(\mathfrak{X}_n)$, $n \in \Lambda$, by restricting $f$ to $X_n$. We can see that $M_{k,\mathfrak{X}}(f)(x) = \max \{M_{k,\mathfrak{X}_n}(f_n)(x), \|f\|_1 / \mu(X)\}$ for $x \in X_n$ and hence, applying H\"older's inequality, we get
\begin{displaymath}
\|M_{k,\mathfrak{X}}(f)\|_p^p \leq \sum_{n \in \Lambda} \|M_{k,\mathfrak{X}_n}(f_n)\|_p^p + \|f\|_1^p \cdot \mu(X)^{1-p} \leq \sum_{n \in \Lambda} \tilde{c}^p \|f_n\|_p^p + \|f\|_p^p = (\tilde{c}^p + 1) \|f\|_p^p. \raggedright \hfill \qed
\end{displaymath}	
\end{prof*}

Let us note here that whenever we want to apply Proposition 1 in this paper, we omit the details related to the proper indexing of the component spaces. We do not even specify $\Lambda$. The only important thing is that each time we use at most countably many spaces. Finally, notice that in the previous related articles, \cite{Ko} and \cite{Ko2}, all the investigated spaces consisted of infinitely many distant parts, say branches, and it was necessary to properly argue that the interactions between the different parts are small enough. Now we can first take a single branch space and study the behavior of the associated maximal operators, and then, by using Proposition 1, combine the branches to build the expected space. Such a strategy seems more natural and simplifies calculations. This will be particularly evident in Section 3, where the so-called basic spaces will be introduced.
	
\section{Results for single $k$}	
	
For a fixed metric measure space $\mathfrak{X}$ and $k \in [1, \infty)$ denote by $P_{k,\rm s}^{\rm c}$ and $P_{k,\rm s}$ the sets consisting of such $p \in [1, \infty]$ for which the associated operators, $M^{\rm c}_k$ and $M_k$, are of strong type $(p,p)$, respectively. Similarly, let $P_{k,\rm w}^{\rm c}$ and $P_{k,\rm w}$ consist of such $p \in [1, \infty]$ for which $M^{\rm c}_k$ and $M_k$ are of weak type $(p,p)$, respectively. Then\smallskip

\begin{enumerate}[label=(\roman*)]
	\item each of the four sets is of the form $\{\infty\}$, $[p_0, \infty]$ or $(p_0,\infty]$, for some $p_0 \in [1, \infty)$ (this is a natural consequence of the $L^\infty$ boundedness of the considered operators and the Marcinkiewicz interpolation theorem);\smallskip
	
	\item we have the following inclusions
	\begin{displaymath}
	P_{k,\rm s} \subset P_{k,\rm s}^{\rm c}, \quad P_{k,\rm w} \subset P_{k,\rm w}^{\rm c}, \quad P_{k,\rm s}^{\rm c} \subset P_{k,\rm w}^{\rm c} \subset \overline{P_{k,\rm s}^{\rm c}}, \quad P_{k,\rm s} \subset P_{k,\rm w} \subset \overline{P_{k,\rm s}},
	\end{displaymath}
	where $\overline{E}$ denotes the closure of $E$ in the usual topology of $\mathbb{R} \cup \{ \infty \}$;
	
	\item if $k \geq 2$, then $P_{k,\rm w}^{\rm c} = [1, \infty]$ (see \cite{Sa}, \cite{Te} and \cite{St1});
	
	\item if $k \geq 3$, then $P_{k,\rm w} = [1, \infty]$ (see \cite{NTV} and \cite{St1}).
\end{enumerate}

In this section we study the possible configurations of the sets $P_{k,\rm s}^{\rm c}$, $P_{k,\rm s}$, $P_{k,\rm w}^{\rm c}$ and $P_{k,\rm w}$ for a fixed $k \in [1, \infty)$. Let us mention at this point that sometimes, if it is significant, we indicate the space with respect to which the given set was built (for example, we write $P_{k,\rm s}^{\rm c}(\mathfrak{X})$ instead of $P_{k,\rm s}^{\rm c}$). Moreover, if $k = 1$, then we write shortly $P_{\rm s}^{\rm c}$ instead of $P_{k,\rm s}^{\rm c}$ and so on. It is worth noting here that the case $k=1$ has been completely investigated in \cite{Ko} (see also \cite{Ko2}, where the restricted weak type inequalities was taken into account). Now we will do a similar analysis for each $k \geq 1$. Namely, we will prove the following.

\begin{theorem}
	Fix $k \in [1, \infty)$. Let $P_{k,\rm s}^{\rm c}$, $P_{k,\rm s}$, $P_{k,\rm w}^{\rm c}$ and $P_{k,\rm w}$ be such that conditions $\rm (i)$ and $\rm (ii)$ (and $\rm (iii)$ or $\rm (iv)$, if necessary) hold. Then there exists a (non-doubling) metric measure space for which the associated modified Hardy--Littlewood maximal operators, centered $M^{\rm c}_k$ and non-centered $M_k$, satisfy\smallskip
	\begin{itemize}
		\item $M^{\rm c}_k$ is of strong type $(p,p)$ if and only if $p \in P_{k,\rm s}^{\rm c}$,\smallskip
		\item $M_k$ is of strong type $(p,p)$ if and only if $p \in P_{k,\rm s}$,\smallskip
		\item $M^{\rm c}_k$ is of weak type $(p,p)$ if and only if $p \in P_{k,\rm w}^{\rm c}$,\smallskip
		\item $M_k$ is of weak type $(p,p)$ if and only if $p \in P_{k,\rm w}$.\smallskip
	\end{itemize} 
\end{theorem}
We will prove Theorem 1 in Section 2.3. To do that we need a few auxiliary lemmas, which will be formulated in Sections 2.1 and 2.2. Let us also comment that the most interesting case concerns $k \in [1,3)$. If $k \geq 3$, then we have only three possibilities depending on whether $M_k^{\rm c}$ and $M_k$ are of strong type $(1,1)$ or not.

\subsection{First and second generation spaces} To prove Theorem 1 we use results obtained in \cite{Ko}, where some specific non-doubling metric measure spaces, so called first and second generation spaces, were considered. We give only a brief description of these spaces and do not go far into details, kindly asking the reader to consult \cite{Ko} if necessary.

Now we present the construction process for first generation spaces. 
Let $\tau = (\tau_n)_{n \in \mathbb{N}}$ be a fixed sequence of positive integers. Define 
\begin{displaymath}
	X_{\tau} = \{x_n \colon n \in \mathbb{N}\} \cup \{x_{ni} \colon  i=1, \dots, \tau_n, n \in \mathbb{N}\},
\end{displaymath}
where all elements $x_n, x_{ni}$ are pairwise different. We introduce the metric $\rho = \rho_\tau$ determining the distance between two different elements $x$ and $y$ by the formula
\begin{displaymath}
	\rho(x,y) = \left\{ \begin{array}{rl}
		1 & \textrm{if } x_n \in \{x,y\} \subset S_n \textrm{ for some } n \in \mathbb{N},  \\
		2 & \textrm{otherwise,} \end{array} \right. 
\end{displaymath}
where $S_n = \{ x_n, x_{n1}, \dots , x_{n\tau_n} \}$. 
Finally, we define the measure $\mu = \mu_{\tau, F}$ on $X_{\tau}$ by letting $\mu(\{x_n\}) = d_n$ and $\mu(\{x_{ni}\}) = d_nF(n,i)$, where $F > 0$ is a given function and $d = (d_n)_{n \in \mathbb{N}}$ is an appropriate sequence of strictly positive numbers with $d_1 = 1$ and $d_n$ chosen (uniquely!) in such a way that $\mu(S_n) = \mu(S_{n-1})/2$, $n \geq 2$. Note that this implies $\mu(X_\tau) < \infty$. Moreover, observe that $\mu$ is non-doubling.

Next, we describe second generation spaces. Let $\tau^* = (\tau^*_n)_{n \in \mathbb{N}}$ be a fixed sequence of positive integers. Define
\begin{displaymath}
Y_{\tau^*} = \{y_n \colon n \in \mathbb{N}\} \cup \{y_{ni}, y_{ni}' \colon i = 1, \dots, \tau^*_n, n \in \mathbb{N} \},
\end{displaymath}
where all elements $y_n, y_{ni}, y_{ni}'$ are pairwise different. We introduce the metric $\rho = \rho_{\tau^*}$ determining the distance between two different elements $x$ and $y$ by the formula
\begin{displaymath}
\rho(x,y) = \left\{ \begin{array}{rl}
1 & \textrm{if } \{x,y\} = T_{ni} \textrm{ or } y_n \in \{x,y\} \subset T_n \setminus T_n ' \textrm{ for some } n \in \mathbb{N}, \ i \in \{1, \dots, \tau^*_n\}, \\
2 & \textrm{otherwise,} \end{array} \right. 
\end{displaymath}
where $T_n = \{ y_n, y_{n1}, \dots , y_{n\tau^*_n}, y_{n1}', \dots , y_{n\tau^*_n}' \}$, $T_n ' = \{y_{n1}', \dots , y_{n\tau^*_n}'\}$ and $T_{ni}=\{y_{ni}, y_{ni}'\}$. Finally, we define the measure $\mu = \mu_{\tau^*, F^*}$ by letting $\mu(\{y_n\}) = d^*_n$, $\mu(\{y_{ni}\}) = d^*_n / \tau^*_n$ and $\mu(\{y_{ni}'\}) = d^*_n
F^*(n,i)$, where $F^* > 0$ is a given function and $d^* = (d^*_n)_{n \in \mathbb{N}}$ is an appropriate sequence of strictly positive numbers with $d^*_1 = 1$ and $d^*_n$ chosen (uniquely!) in such a way that $\mu(T_n) = \mu(T_{n-1})/2$, $n \geq 2$. Note that this implies $\mu(Y_{\tau^*}) < \infty$ and observe that $\mu$ is non-doubling. In addition, as it is proven in \cite{Ko}, for each second generation space the associated centered maximal operator is of strong type $(1,1)$.

In \cite{Ko} described are all possible configurations of the sets $P_{\rm s}^{\rm c}$, $P_{\rm s}$, $P_{\rm w}^{\rm c}$ and $P_{\rm w}$, by using the first and second generation spaces and some mixed spaces, which are constructed in the spirit of Proposition 1 (in this process we combine two component spaces and the distance between elements belonging to different pieces equals $2$). Note that for any such a space $\mathfrak{X}$ the metric $\rho$ takes only two non-zero values, $1$ and $2$. Therefore, in this case, for any $k \in [1,2)$ the operators $M_k^{\rm c}$ and $M_k$ are identical with $M^{\rm c}$ and $M$, respectively. The key point here is that for $k \in [1,2)$ we can find $r>1$ such that $kr \leq 2$. As a result, we obtain the equalities $P_{k,\rm s}^{\rm c}(\mathfrak{X})=P_{\rm s}^{\rm c}(\mathfrak{X})$, $P_{k,\rm s}(\mathfrak{X})=P_{\rm s}(\mathfrak{X})$, $P_{k,\rm w}^{\rm c}(\mathfrak{X})= P_{\rm w}^{\rm c}(\mathfrak{X})$ and $P_{k,\rm w}(\mathfrak{X})=P_{\rm w}(\mathfrak{X})$. 

The situation when $k \in [2,3)$ is different. Namely, in this case, for any ball $B$ consisting of at least two points the ball $kB$ coincides with the whole space. This fact causes that $M^{\rm c}_k$ and $M_k$ are of strong type $(1,1)$. However, a slight modification of the metric used in the construction of second generation spaces will allow us to obtain more subtle results. 

\begin{lemma}
	Fix $k \in [2, 3)$. Let $\mathfrak{Y} = (Y, \rho, \mu)$ be a second generation space. Define the metric $\rho'$ determining the distance between two different elements $x, y \in Y$ by the formula
	\begin{displaymath}
	\rho'(x,y) = \left\{ \begin{array}{rl}
	1 & \textrm{if } \rho(x,y)=1,  \\
	2 & \textrm{if there exists } z \in Y \textrm{ such that } \rho(x,z)= \rho(y,z)=1, \\
	3 & \textrm{otherwise.} \end{array} \right. 
	\end{displaymath}
	Then for the space $\mathfrak{Y}' = (Y, \rho', \mu)$ we have $P_{k,\rm s}^{\rm c}(\mathfrak{Y}') = P_{k,\rm w}^{\rm c}(\mathfrak{Y}') = [1, \infty]$, while $P_{k,\rm s}(\mathfrak{Y}') = P_{\rm s}(\mathfrak{Y})$ and $P_{k,\rm w}(\mathfrak{Y}') = P_{\rm w}(\mathfrak{Y})$.
\end{lemma}

\begin{prof*}

First of all, let us emphasize that $\rho'$ is well-defined. Indeed, it can be easily seen that there is no set $\{x, y, z\} \subset Y$ satisfying
\begin{displaymath}
\rho(x,y) =  \rho(x, z) = \rho(y,z) = 1,
\end{displaymath} 
and thus the first two conditions in the definition of $\rho'$ cannot happen at the same time.

Now observe that $L^1(\mathfrak{Y})$ and $L^1(\mathfrak{Y}')$ are equal as Banach spaces. Moreover, for any $f \in L^1(\mathfrak{Y})$ we have $M_{k,\mathfrak{Y}'}^{\rm c}(f) \leq M^{\rm c}_\mathfrak{Y}(f)$ and $M_{k,\mathfrak{Y}'}(f) \leq M_\mathfrak{Y}(f)$. Indeed, in the case of the centered operators, suppose that $f \geq 0$ and fix $y \in Y$. If $r \leq 2$, then we have $B_{\rho'}(y,kr) \supset B_{\rho'}(y,r) = B_{\rho}(y,r)$, which implies
\begin{displaymath}
\frac{1}{\mu(B_{\rho'}(y,kr))} \sum_{x \in B_{\rho'}(y,r)} f(x) \mu(\{x\})  \leq \frac{1}{\mu(B_{\rho}(y,r))} \sum_{x \in B_{\rho}(y,r)} f(x) \mu(\{x\}) \leq M^{\rm c}_\mathfrak{Y}(f)(y).
\end{displaymath}
On the other hand, if $r > 2$, then we have $B_{\rho'}(y,kr) = Y$, which implies
\begin{displaymath}
\frac{1}{\mu(B_{\rho'}(y,kr))} \sum_{x \in B_{\rho'}(y,r)} f(x) \mu(\{x\})  \leq \frac{1}{\mu(Y)} \sum_{y \in Y} f(y) \mu(\{y\}) \leq M^{\rm c}_\mathfrak{Y}(f)(y).
\end{displaymath}
This gives $M_{k,\mathfrak{Y}'}^{\rm c}(f) \leq M^{\rm c}_\mathfrak{Y}(f)$ and the second claimed estimate is verified analogously. Hence, we obtain the following equalities and inclusions
\begin{displaymath}
P_{k,\rm s}^{\rm c}(\mathfrak{Y}') = P_{k,\rm w}^{\rm c}(\mathfrak{Y}') = [1, \infty], \quad P_{k,\rm s}( \mathfrak{Y}') \supset P_{\rm s}(\mathfrak{Y}), \quad P_{k,\rm w}(\mathfrak{Y}') \supset P_{\rm w}(\mathfrak{Y}).
\end{displaymath}
Let us point out here that, in particular, we have $P_{\rm s}^{\rm c}(\mathfrak{Y}') = [1, \infty]$, which is quite peculiar since in most typical settings $M^{\rm c}$ is not of strong type $(1,1)$.
 
Now, it remains to show that if $M_\mathfrak{Y}$ is not of strong (respectively weak) type $(p,p)$ for some $p \geq 1$, then $M_{k, \mathfrak{Y}'}$ fails to be of strong (respectively weak) type $(p,p)$. Our strategy is as follows. We present quickly the argument that was used in \cite{Ko} to obtain certain property of $M_\mathfrak{Y}$ and then try to convince the reader that the situation is very similar in the context of $M_{\mathfrak{Y}'}$ instead. For clarity we describe only the case related to the strong type $(p,p)$ inequalities. 

Recall that each time when it was shown that the non-centered operator associated with the second generation space $\mathfrak{Y}$ is not of strong type $(p,p)$, the functions $f_n = \delta_{y_n}$, $n \in \mathbb{N}$, were considered. Then, the functions $M_\mathfrak{Y} (f_n)$ were estimated from below by$\colon$
\begin{itemize}
\item the average value of $f_n$ on the ball centered at $y_{ni}$ with the radius $3/2$ (denoted by $A_{B_\rho(y_{ni},3/2)}f_n$) for the points $y_{ni}'$, $i = 1, \dots, \tau_n$,
\item 0 for the other points,	
\end{itemize}   
and finally it turned out that 
\begin{displaymath}
\lim_{n \rightarrow \infty} \frac{\|M_\mathfrak{Y}(f_n)\|_p^p}{\|f_n\|_p^p} \geq \lim_{n \rightarrow \infty} \frac{ \sum_{i=1}^{\tau_n} (A_{B_\rho(y_{ni},3/2)}f_n)^p \mu(\{y_{ni}' \})
	}{\|f_n\|_p^p} = \infty.
\end{displaymath}

Let us assume that the above estimate holds for some $p \in [1, \infty)$. Take $r>1$ such that $kr \leq 3$ and see that $B_{\rho'}(y_{ni},r) = B_\rho(y_{ni},3/2)$ and 
\begin{displaymath}
\mu(B_{\rho'}(y_{ni},kr)) = \mu\big(\{y_{ni}'\} \cup \big(T_n \setminus T_n'\big) \big) \leq 2 \mu(B_\rho(y_{ni},3/2)).
\end{displaymath}
 This implies $M_{k,\mathfrak{Y}'}(f_n)(y_{ni}') \geq \frac{1}{2} A_{B_\rho(y_{ni},3/2)}f_n$ and hence
\begin{displaymath}
\lim_{n \rightarrow \infty} \frac{\|M_{k,\mathfrak{Y}'}(f_n)\|_p}{\|f_n\|_p} = \infty. \eqno \qed
\end{displaymath} 
\end{prof*}

\subsection{Endpoint cases} Now, we turn our attention to certain specific situations in which the operators $M^{\rm c}_k$ or $M_k$ are not of strong type $(1,1)$ for some $k \geq 2$ or $k \geq 3$, respectively. We begin with a construction of some auxiliary metric measure spaces called by us the $\textit{segment-type spaces}$ and then we prove two lemmas related to them.

Let $d = \{d_{n,i}  \colon i=1, \dots, n ,\ n \in \mathbb{N}\}$ be a fixed triangular matrix of strictly positive numbers such that $\sum_{i=1}^n d_{n,i} \leq 1$, $n \in \mathbb{N}$. Define 
\begin{displaymath}
X = \{x_{n,i} \colon i=0, \dots, n, \ n \in \mathbb{N}\},
\end{displaymath}
where all elements $x_{n,i}$ are pairwise different (and located on the plane, say). By $S_n$ we denote the branch $S_n = \{ x_{n,0}, x_{n,1}, \dots , x_{n,n} \}$. Thus, $X = \bigcup_{n=1}^\infty S_n$ is the disjoint union of the family of branches. We define the metric $\rho = \rho_d$ determining the distance between two different elements $x, y \in X$ by the formula
\begin{displaymath}
\rho(x,y) = \left\{ \begin{array}{rl}
\sum_{i=j+1}^k d_{n,i} & \textrm{if } \{x,y\} = \{x_{n,j}, x_{n,k}\} \textrm{ for some } 0 \leq j < k \leq n, n \in \mathbb{N},  \\
1 & \textrm{otherwise.} \end{array} \right. 
\end{displaymath}
Observe that $\rho$ is determined uniquely by $d$ and clearly $\rm{diam}$$(X) = 1$ holds from the definition. Figure 1 shows a model of the space $(X, \rho)$.

\begin{figure}[H]
	\begin{tikzpicture}
	[scale=.8,auto=left,every node/.style={circle,fill,inner sep=2pt}]
	
	\node[label={[yshift=-1cm]$x_{1,0}$}] (m0) at (2,8) {};
	\node[label={[yshift=-1cm]$x_{1,1}$}] (m1) at (14,8) {};
	
	\node[label={[yshift=-1cm]$x_{2,0}$}] (n0) at (2,6) {};
	\node[label={[yshift=-1cm]$x_{2,1}$}] (n1) at (8,6) {};
	\node[label={[yshift=-1cm]$x_{2,2}$}] (n2) at (14,6) {};
	
	\node[label={[yshift=-1cm]$x_{3,0}$}] (o0) at (2,4) {};
	\node[label={[yshift=-1cm]$x_{3,1}$}] (o1) at (5,4) {};
	\node[label={[yshift=-1cm]$x_{3,2}$}] (o2) at (8,4) {};
	\node[label={[yshift=-1cm]$x_{3,3}$}] (o3) at (14,4) {};
	
	\node[label={[yshift=-1cm]$x_{4,0}$}] (p0) at (2,2) {};
	\node[label={[yshift=-1cm]$x_{4,1}$}] (p1) at (3.5,2) {};
	\node[label={[yshift=-1cm]$x_{4,2}$}] (p2) at (5,2) {};
	\node[label={[yshift=-1cm]$x_{4,3}$}] (p3) at (8,2) {};
	\node[label={[yshift=-1cm]$x_{4,4}$}] (p4) at (14,2) {};
	
	
	\draw (m0) -- (m1) node [midway, fill=white, above=-5pt] {$d_{1,1}$};
	
	\draw (n0) -- (n1) node [midway, fill=white, above=-5pt] {$d_{2,1}$};
	\draw (n1) -- (n2) node [midway, fill=white, above=-5pt] {$d_{2,2}$};
	
	\draw (o0) -- (o1) node [midway, fill=white, above=-5pt] {$d_{3,1}$};
	\draw (o1) -- (o2) node [midway, fill=white, above=-5pt] {$d_{3,2}$};
	\draw (o2) -- (o3) node [midway, fill=white, above=-5pt] {$d_{3,3}$};
	
	\draw (p0) -- (p1) node [midway, fill=white, above=-5pt] {$d_{4,1}$};
	\draw (p1) -- (p2) node [midway, fill=white, above=-5pt] {$d_{4,2}$};
	\draw (p2) -- (p3) node [midway, fill=white, above=-5pt] {$d_{4,3}$};
	\draw (p3) -- (p4) node [midway, fill=white, above=-5pt] {$d_{4,4}$};
	
	\foreach \from/\to in {m0/m1, n0/n2, o0/o3, p0/p4}
	\draw (\from) -- (\to);
	
	\draw [thick, dotted] (2,0.8) -- (2,0.5);
	\draw [thick, dotted] (14,0.8) -- (14,0.5);
	\draw [thick, dotted] (7.85,0.6) -- (8.15,0.6);
	\end{tikzpicture}
	\caption{The model of the space $(X, \rho)$.}
\end{figure}
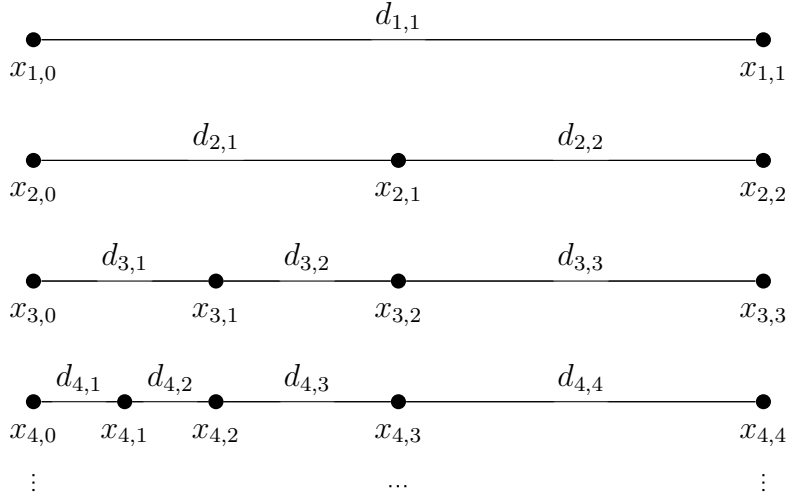

We define the measure $\mu = \mu_{F}$ on $X$ by letting $\mu(\{x_{n,i}\}) = F(n,i)$, where $F>0$ is a given function satisfying $ \sum_{i=0}^n F(n,i) \leq 2^{-n}$, $n \in \mathbb{N}$ (this implies $\mu(X) < \infty$). Observe that $(X, \rho, \mu)$ is non-doubling. Indeed, fix $\epsilon > 0$ and let $n_0 = n_0(\epsilon)$ be such that $\mu(S_{n_0}) < \epsilon$. Then we have $B(x_{n_0,0},2/3) \subset S_{n_0}$ which implies $\mu(B(x_{n_0,0},2/3)) < \epsilon$, while $\mu(B(x_{n_0,0},4/3)) = \mu(X)$.

From now on we shall use the sign $|E|$ instead of $\mu(E)$ for $E \subset X$. It will be clear from the context when the symbol $|\cdot|$ refers to the measure and when it denotes the absolute value sign. 

\begin{lemma}
	Fix $k \geq 2$ and let $\mathfrak{X} = (X, \rho, \mu)$ be the segment-type space with $d_{n,i} = (k+1)^{i-n-1}$ and $F(n,i) = 2^{-n} (n+1)^{-1}$, $i=1, \dots n$, $n \in \mathbb{N}$. Then $P_{k,\rm s}^{\rm c}(\mathfrak{X}) = P_{k,\rm s}(\mathfrak{X}) = (1, \infty]$ and $P_{k,\rm w}^{\rm c}(\mathfrak{X}) = P_{k,\rm w}(\mathfrak{X}) = [1, \infty]$.
\end{lemma}

\begin{prof*}
	At the beginning, note that it suffices to show that $M_k^{\rm c}$ is not of strong type $(1,1)$, while $M_k$ is of weak type $(1,1)$. Observe that for $j= 1, \dots, n-1$, $n \in \mathbb{N}$, we have
	\begin{displaymath}
	\sum_{i=1}^{j} d_{n,i} = \sum_{i=1}^{j} \frac{1}{(k+1)^{n-i+1}} < \sum_{i=1}^{\infty} \frac{1}{(k+1)^{n-j+i}} = \frac{1}{k(k+1)^{n-j}} = \frac{d_{n,j+1}}{k} \leq \frac{1}{k}.
	\end{displaymath}
First we show that $M_k^{\rm c}$ is not of strong type $(1,1)$. Let $f_n = \delta_{x_{n,0}}$, $n \geq 1$. Then $\|f_n\|_1 = |\{x_{n,0}\}|$. For any $j= 1, \dots, n-1$ we can find $r=r(j)$ such that $B(x_{n,j},r) = B(x_{n,j}, kr) = \{x_{n,0}, x_{n,1}, \dots, x_{n,j}\}$ and hence $M_k^{\rm c}f_n(x_{n,j}) \geq \frac{1}{j+1}$ for that $j$. This implies
	\begin{displaymath}
	\limsup_{n \rightarrow \infty} \frac{\|M_k^{\rm c}f_n\|_1}{\|f_n\|_1} \geq \limsup_{n \rightarrow \infty} \frac{\sum_{j=1}^{n-1}M_k^{\rm c}f_n(x_{n,j})|\{x_{n,j}\}|}{|\{x_{n,0}\}|} \geq \lim_{n \rightarrow \infty} \sum_{j=1}^{n-1} \frac{1}{j+1}= \infty.
	\end{displaymath} 
	
Now, it remains to show that $M_k$ is of weak type $(1,1)$. Let $f \in L^{1}(\mathfrak{X})$, $f \geq 0$, and consider $\lambda > 0$ such that $E_\lambda (M_kf) \neq \emptyset$. If $\lambda < \|f\|_1 / |X|$, then $\lambda |E_\lambda(M_kf)| / \|f\|_1 \leq 1$ follows. Therefore, from now on assume that $\lambda \geq \|f\|_1 / |X|$. With this assumption, if for some $x \in S_n$ we have $M_kf(x) > \lambda$, then any ball $B$ containing $x$ and realizing $\sum_{x \in B} f(x) |\{x\}| / |kB| > \lambda$ must be a subset of $S_n$. Moreover, because of the linear structure of $S_n$, any ball $B \subset S_n$ is of the form $B = \{ x_{n,i}, x_{n,i+1}, \dots, x_{n,j} \}$ for some $0 \leq i \leq j \leq n$. Take any $n \in \mathbb{N}$ such that $E_\lambda(M_kf) \cap S_n \neq \emptyset$ and consider $\mathcal{B} = \mathcal{B}(n) = \{B \subset S_n \colon \sum_{x \in B} f(x) |\{x\}| / |kB| > \lambda\}$ which forms a cover of $E_\lambda(M_kf) \cap S_n$. By using the fact that each element of $\mathcal{B}$ has the form described above we can find a subcover $\mathcal{B}'$ such that each $x \in E_\lambda(M_kf) \cap S_n$ belongs to at most two elements of $\mathcal{B}'$. Therefore
\begin{align*}
\lambda |E_\lambda(M_kf) \cap S_n| & \leq \sum_{B \in \mathcal{B}'} \lambda |B| \leq \sum_{B \in \mathcal{B}'} \Big( \sum_{x \in B} f(x) |\{x\}| / |kB| \Big) |B| \\
& \leq \sum_{B \in \mathcal{B}'} \sum_{x \in B} f(x) |\{x\}| \leq 2 \sum_{x \in E_\lambda(M_kf) \cap S_n} f(x) |\{x\}|,
\end{align*}
and, consequently, $\lambda |E_\lambda(M_kf)| \leq 2 \|f\|_1. \raggedright \hfill \qed$ 	
\end{prof*}

\begin{lemma}
	Fix $k \geq 3$ and let $\mathfrak{X} = (X, \rho, \mu)$ be the segment-type space with $d_{n,i} = (k-1/2)^{i-n-1}$, $i=1, \dots n$, $n \in \mathbb{N}$, and $F(n,i)$ chosen (uniquely) in such a way that $F(n,n) = 2^{-n-1}$ and $F(n,i) = F(n,i+1) / 2^{i+1}$ for $i \in \{0, \dots, n-1\}$, $n \in \mathbb{N}$.
	Then $P_{k,\rm s}(\mathfrak{X}) = (1, \infty]$ and $P_{k,\rm s}^{\rm c}(\mathfrak{X}) = P_{k,\rm w}^{\rm c}(\mathfrak{X}) = P_{k,\rm w}(\mathfrak{X}) = [1, \infty]$.
\end{lemma}

\begin{prof*}
	Note that, since $k \geq 3$, $M_k$ is of weak type $(1,1)$. Hence, it suffices to show that $M_k$ is not of strong type $(1,1)$, while $M_k^{\rm c}$ is of strong type $(1,1)$. Observe that $\sum_{i=1}^{n} d_{n,i} < 1$ and $\sum_{i=1}^{j} d_{n,i} < d_{n,j+1}$, $j=0, 1, \dots, n-1$. Moreover, one can see that $\sum_{i=1}^{n} F(n,i) < 2^{-n}$ and $\sum_{i=1}^{j} F(n,i) < F(n,j+1)$, $j=0, 1, \dots, n-1$.
	
	First we show that $M_k$ is not of strong type $(1,1)$. Let $f_n = \delta_{x_{n,0}}$, $n \geq 1$. Then $\|f_n\|_1 = |\{x_{n,0}\}|$. Since $\sum_{i=1}^{j-1} d_{n,i} < d_{n,j} < d_{n,j+1} / (k-1)$ for any $j= 1, \dots, n-1$, then we can find $r=r(j)$ such that $B(x_{n,j-1},r) = B(x_{n,j-1}, kr) = \{x_{n,0}, x_{n,1}, \dots, x_{n,j}\}$ and hence $M_kf_n(x_{n,j}) \geq \frac{|\{x_{n,0}\}|}{2|\{x_{n,j}\}|}$ for that $j$. This implies
	\begin{displaymath}
	\limsup_{n \rightarrow \infty} \frac{\|M_kf_n\|_1}{\|f_n\|_1} \geq \limsup_{n \rightarrow \infty} \frac{\sum_{j=1}^{n-1}M_kf_n(x_{n,j})|\{x_{n,j}\}|}{|\{x_{n,0}\}|} \geq \lim_{n \rightarrow \infty} \frac{n-1}{2}= \infty.
	\end{displaymath} 
	Now, it remains to show that $M^{\rm c}_k$ is of strong type $(1,1)$. Let $f \in L^{1}(\mathfrak{X})$, $f \geq 0$. Observe that $d_{n,j} \, k > d_{n,j+1}$ and hence for any ball $B$ centered at $x_{n,j}$, $1 \leq j \leq n-1$, if $x_{n,j-1} \in B$, then $x_{n,j+1} \in kB$. Therefore we have the estimate
	\begin{displaymath}
	M^{\rm c}_kf(x_{n,j}) \leq f(x_{n,j}) + \frac{ \sum_{i=1}^{j} f(x_{n,i}) |\{x_{n,i} \}|}{|\{x_{n,j+1}\}|} + \max\{f(x_{n,j+1}), \dots, f(x_{n,n})\} + \frac{\|f\|_1}{|X|},
	\end{displaymath}
	which implies
	\begin{align*}
	\sum_{n=1}^{\infty} \sum_{j=0}^{n} M^{\rm c}_kf(x_{n,j}) |\{x_{n,j} \}| & \leq \sum_{n=1}^{\infty} \sum_{j=0}^{n} f(x_{n,j}) |\{x_{n,j} \}| \Big( 2 + \sum_{i=j}^{n} \frac{|\{x_{n,i} \}|}{|\{x_{n,i+1} \}|} + \sum_{i=0}^{j-1} \frac{|\{x_{n,i} \}|}{|\{x_{n,j} \}|} \Big) \\
	& \leq \sum_{n=1}^{\infty} \sum_{j=0}^{n} f(x_{n,j}) |\{x_{n,j} \}| \Big( 2 + 1 + 1 \Big) = 4 \, \|f\|_1. \qquad \qquad \qquad \qed
	\end{align*}
\end{prof*}
\subsection{Proof of Theorem 1}

	At the beginning note that if $P_{k,\rm s}^{\rm c} = P_{k,\rm s} = P_{k,\rm w}^{\rm c} = P_{k,\rm w} = [1, \infty]$, then we can find a first generation space $\mathfrak{X}$ for which $P_{\rm s}^{\rm c}(\mathfrak{X}) = P_{\rm s}(\mathfrak{X}) = P_{\rm w}^{\rm c}(\mathfrak{X}) = P_{\rm w}(\mathfrak{X}) = [1, \infty]$, and hence we also have $P_{k,\rm s}^{\rm c}(\mathfrak{X}) = P_{k,\rm s}(\mathfrak{X}) = P_{k,\rm w}^{\rm c}(\mathfrak{X}) = P_{k,\rm w}(\mathfrak{X}) = [1, \infty]$ for every $k \geq 1$. Therefore, from now on, assume that $P_{k,\rm s}$ (and possibly the other sets) is a proper subset of $[1, \infty]$. We shall consider the cases$\colon$ $k \in [1,2)$, $k \in [2,3)$ and $k \geq 3$.
	
	First, suppose that $k \in [1,2)$. Then we assume that the sets $P_{k,\rm s}^{\rm c}$, $P_{k,\rm s}$, $P_{k,\rm w}^{\rm c}$ and $P_{k,\rm w}$ satisfy $\rm (i)$ and $\rm (ii)$. We can find a (first or second generation, or mixed) space $\mathfrak{X}$ for which $P_{\rm s}^{\rm c}( \mathfrak{X}) = P_{k,\rm s}^{\rm c}$, $P_{\rm s}(\mathfrak{X}) = P_{k,\rm s}$, $P_{\rm w}^{\rm c}(\mathfrak{X}) = P_{k,\rm w}^{\rm c}$ and $P_{\rm w}(\mathfrak{X}) = P_{k,\rm w}$, and using the observation from Section 2.1 we can see that the same equalities with $k$ instead of $1$ hold.
	
	Next, suppose that $k \in [2,3)$. Then we assume that the sets $P_{k,\rm s}^{\rm c}$, $P_{k,\rm s}$, $P_{k,\rm w}^{\rm c}$ and $P_{k,\rm w}$ satisfy $\rm (i)$, $\rm (ii)$, and $\rm (iii)$. We can find a second generation space $\mathfrak{Y}$ for which $P_{\rm s}^{\rm c}(\mathfrak{Y}) = P_{\rm w}^{\rm c}(\mathfrak{Y}) = [1, \infty] = P_{k,\rm w}^{\rm c}$, $P_{\rm s}(\mathfrak{Y}) = P_{k,\rm s}$ and $P_{\rm w}( \mathfrak{Y}) = P_{k,\rm w}$, and therefore we obtain the adequate equalities with $P_{\rm s}^{\rm c}(\mathfrak{Y})$, $P_{\rm s}(\mathfrak{Y})$, $P_{\rm w}^{\rm c}(\mathfrak{Y})$ and $P_{\rm w}(\mathfrak{Y})$ replacing by $P_{k,\rm s}^{\rm c}(\mathfrak{Y}')$, $P_{k,\rm s}(\mathfrak{Y}')$, $P_{k,\rm w}^{\rm c}(\mathfrak{Y}')$ and $P_{k,\rm w}(\mathfrak{Y}')$, respectively, where $\mathfrak{Y}'$ is the modification of $\mathfrak{Y}$ considered in Lemma 1. If $P_{k,\rm s}^{\rm c} = [1, \infty]$, then the expected space may be chosen to be just $\mathfrak{Y}'$. In the other case (that is, if $P_{k,\rm s}^{\rm c} = (1, \infty]$) we use Proposition 1 with $k_0=k$ to combine $\mathfrak{Y}'$ with the appropriate segment-type space considered in Lemma 2. After that, we obtain the new space, say $\mathfrak{Z}$, such that $P_{k,\rm s}^{\rm c}(\mathfrak{Z}) = P_{k,\rm s}^{\rm c}$, $P_{k,\rm s}(\mathfrak{Z}) = P_{k,\rm s}$, $P_{k,\rm w}^{\rm c}(\mathfrak{Z}) = P_{k,\rm w}^{\rm c}$ and $P_{k,\rm w}(\mathfrak{Z}) = P_{k,\rm w}$.
	
	Finally, suppose that $k \geq 3$ and the sets $P_{k,\rm s}^{\rm c}$, $P_{k,\rm s}$, $P_{k,\rm w}^{\rm c}$ and $P_{k,\rm w}$ satisfy $\rm (i)$, $\rm (ii)$, $\rm (iii)$ and $\rm (iv)$. If $P_{k,\rm s}^{\rm c} = P_{k,\rm s} = (1, \infty]$, then the expected space may be chosen to be the suitable space considered in Lemma 2. In the other case, if $P_{k,\rm s}^{\rm c} = [1, \infty]$ and $P_{k,\rm s} = (1, \infty]$, then the expected space may be chosen to be the suitable space considered in Lemma 3. Thus, the proof of Theorem 1 is complete. 

\section{Results for varying $k$}

For a fixed metric measure space $\mathfrak{X}$ and parameters $p \in [1, \infty]$ and $k \geq 1$ we denote by $c(k,p) = c_\mathfrak{X}(k,p)$ the best constant in the weak type $(p,p)$ inequality for the associated maximal operator $M_k$ (if $M_k$ is not of weak type $(p,p)$, then we write $c(k,p)= \infty$). Similarly, we define $c^{\rm c}(k,p)$ with $M_k^{\rm c}$ replacing $M_k$. In this section we try to study the behavior of these functions, in particular with regard to when they are finite or not. With this in mind, let us define auxiliary functions $h^{\rm c}(k) = \inf\{p \colon c^{\rm c}(k,p) < \infty\}$ and $h(k) = \inf\{p \colon c(k,p) < \infty\}$. Since $M_2^{\rm c}$ and $M_3$ are of weak type $(1,1)$ we can assume that the domains of $h^{\rm c}$ and $h$ are $[1,2]$ and $[1,3]$, respectively. We have the following properties
\begin{enumerate}[label=(\alph*)]
	\item $h^{\rm c} \colon [1,2] \rightarrow [1, \infty]$ and $h \colon [1,3] \rightarrow [1, \infty]$,
	\item $h^{\rm c}$ and $h$ are non-increasing,
	\item $h(k) \geq h^{\rm c}(k)$ for $k \in [1,2]$,
	\item $h^{\rm c}(2)=h(3) = 1,$
	\item for a fixed $k \in [1,2]$ the set $P_k^{\rm c}$ coincides with $\{\infty\}$, if $h^{\rm c}(k) = \infty$, and takes the form $(h^{\rm c}(k), \infty]$ or $[h^{\rm c}(k), \infty]$ in the opposite case,
	\item for a fixed $k \in [1,3]$ the set $P_k$ coincides with $\{\infty\}$, if $h(k) = \infty$, and takes the form $(h(k), \infty]$ or $[h(k), \infty]$ in the opposite case.
\end{enumerate}

Our principal motivation is to take arbitrary functions $h^{\rm c}$ and $h$ satisfying $\rm (a)-(d)$ and ask whether it is possible to find a metric measure space $\mathfrak{Z}$ such that $\rm (e)$ and $\rm (f)$ hold for $P_k^{\rm c} = P_k^{\rm c}(\mathfrak{Z})$ and $P_k = P_k(\mathfrak{Z})$, respectively. It turns out that the answer is always positive. However, observe that conditions $\rm (a)-(f)$ do not usually include full information about the finiteness of $c^{\rm c}(k, p)$ and $c(k, p)$. Namely, if we know only the values of $h^{\rm c}$ and $h$, then it is rather impossible to determine whether $c^{\rm c}(k, h^{\rm c}(k))$ and $c(k, h(k))$ are finite or not. Moreover, it seems that, with respect to that, we often have a lot of possible cases and even the characterization of them is a difficult problem which will not be resolved here. Nevertheless, the obtained results may be helpful to find a general principle related to this issue. The main goal in this section is to prove the following.

\begin{theorem}
	Let $h^{\rm c}$ and $h$ be such that conditions $\rm (a)$, $\rm (b)$, $\rm (c)$ and $\rm (d)$ hold. Then there exists a metric measure space $\mathfrak{Z}$ such that the associated modified maximal operators $M_k^{\rm c}$, $k \in [1, 2)$, are of weak type $(p,p)$ if and only if $p > h^{\rm c}(k)$ or $p = \infty$, while the operators $M_k$, $k \in [1, 3)$, are of weak type $(p,p)$ if and only if $p > h(k)$ or $p = \infty$. 
\end{theorem}

In Section 3.1 and Section 3.2 we consider auxiliary structures called $\textit{basic spaces}$. The proof of Theorem 2 is located in Section 3.3. Finally, in Section 3.4 we give examples showing that the situation is quite different if we want to find a space $\mathfrak{Z}$ such that the associated modified maximal operators $M_k^{\rm c}$ and $M_k$ are of weak type $(p,p)$ if and only if $p \geq h^{\rm c}(k)$ and $p \geq h(k)$, respectively. Among other things, there is no counterpart of Theorem 2 with these inequalities replacing $p > h^{\rm c}(k)$ and $p > h(k)$.

\subsection{Basic spaces}
Fix $\tau \in \mathbb{N}$, $1 < d \leq 2$ and $m > 1$. We introduce the basic space $\mathfrak{S} = \mathfrak{S}_{\tau, d, m} = (X, \rho, \mu)$ as follows. Let $X = \{x_0, x_1, \dots, x_\tau\}$. Define $\rho = \rho_d$ by letting $\rho(x,y) = 1$ if $x_0 \in \{x,y\}$ and $\rho(x,y) = d$ otherwise, where $x$ and $y$ are two different elements of $X$. Finally, take $\mu = \mu_m$ satisfying $|\{x_0\}| = 1$ and $|\{x_i\}| = m$ for $i=1, \dots, \tau$. Figure 2 shows a model of the space $\mathfrak{S}$.

\begin{figure}[H]
	\begin{tikzpicture}
	[scale=.7,auto=left,every node/.style={circle,fill,inner sep=2pt}]
	
	\node[label={[yshift=-1cm]$x_0$}] (m0) at (8,1) {};
	\node[label=$x_{1}$] (m1) at (5,3)  {};
	\node[label=$x_{2}$] (m2) at (6.5,3)  {};
	\node[label={[yshift=-0.23cm]$x_{\tau-1}$}] (m3) at (9.5,3)  {};
	\node[label={[yshift=-0.04cm]$x_{\tau}$}] (m4) at (11,3)  {};
	\node[dots, scale=2] (m5) at (8,3)  {...};
	
	\foreach \from/\to in {m0/m1, m0/m2, m0/m3, m0/m4}
	\draw (\from) -- (\to);
	\end{tikzpicture}
	\caption{The model of the space $\mathfrak{S}$.}
\end{figure}
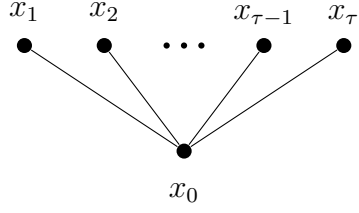
Note that we can explicitly describe any ball:
\begin{displaymath}
B(x_0,r) = \left\{ \begin{array}{rl}
\{x_0\} & \textrm{for } 0 < r \leq 1, \\
X & \textrm{for } 1 < r, \end{array} \right.
\end{displaymath} 
\noindent and, for $i \in \{1, \dots, \tau\}$,
\begin{displaymath}
B(x_i,r) = \left\{ \begin{array}{rl}
\{x_i\} & \textrm{for } 0 < r \leq 1, \\
\{x_0, x_i\} & \textrm{for } 1 < r \leq d,  \\
X & \textrm{for } d < r. \end{array} \right.
\end{displaymath} 

\begin{lemma}
	Let $\mathfrak{S}$ be the metric measure space defined as above. Then
	\begin{displaymath}
	c_\mathfrak{S}(k,p) \simeq c^{\rm c}_\mathfrak{S}(k,p) \simeq \left\{ \begin{array}{rl}
	\max\{1, \tau^{1/p}  m^{1/p-1}\} & \textrm{if } 1 \leq k < d \textrm{ and } p \in [1, \infty),  \\
	1 & \textrm{if } k \geq d \textrm{ or } p = \infty.\end{array} \right. 
	\end{displaymath}
\end{lemma}

\begin{prof*}
	First of all observe that, in view of $c_\mathfrak{S}(k,p) \geq c^{\rm c}_\mathfrak{S}(k,p)$, it suffices to estimate $c_\mathfrak{S}(k,p)$ from above and $c^{\rm c}_\mathfrak{S}(k,p)$ from below by the appropriate terms. Let $f \geq 0$ be a function on $X$. Clearly, $M_k^{\rm c}(f) \geq f$ and hence $c^{\rm c}_\mathfrak{S}(k,p) \geq 1$ for any $k \geq 1$ and $p \in [1, \infty]$. Next, if $k \geq d$ and $p \in [1, \infty)$, then for any ball $B$ consisting of at least two points the ball $kB$ coincides with $X$. Therefore
	\begin{displaymath}
	M_k(f)(x) \leq f(x) + \frac{\|f\|_1}{|X|}, \qquad x \in X.
	\end{displaymath}
	Applying H\"older's inequality we get $\|M_k(f)\|_p^p \leq 2^{p-1} \|f\|_p^p$, $p \in [1, \infty)$, and hence $c_\mathfrak{S}(k,p) \leq 2^{(p-1)/p} \lesssim 1$. Obviously, we also have $c_\mathfrak{S}(k, \infty) \leq 1$ for any $k \geq 1$. From now on assume that $1 \leq k < d$ and $p \in [1, \infty)$. Write $f$ as a sum of the functions $f_1 = f \cdot \chi_{\{x_0\}}$ and $f_2 = f - f_1$. Note that $M_k$ is sublinear, which implies $M_k(f) \leq M_k(f_1) + M_k(f_2)$. We have $M_k(f_1)(x_0) = f_1(x_0)$ and $M_k(f_1)(x_i) \leq  f_1(x_0)/m$ for $i = 1, \dots, \tau$. Then $\|M_k(f_1)\|_p^p \leq (1 + \tau m^{1-p}) \|f_1\|_p^p$. In turn, $M_k(f_2)(x) \leq f_2(x) + \|f_2\|_1 / |X|$, and hence $\|M_k(f_2)\|_p^p \leq 2^{p-1} \|f_2\|_p^p$. Therefore $\|M_k(f)\|_p^p \leq 2^{p-1} (2^{p-1} + 1 + \tau m^{1-p}) \|f\|_p^p$ and we obtain
	\begin{displaymath}
	c_\mathfrak{S}(k,p) \leq \big( 2^{p-1} (2^{p-1} + 1 + \tau m^{1-p}) \big)^{1/p} \lesssim \max\{1, \tau^{1/p}  m^{1/p-1}\}.
	\end{displaymath}
	Finally, if $f = \delta_{\{x_0\}}$, then $\|f\|_p = 1$ and $M^{\rm c}_k(f)(x_i) = 1 / (1+m) \geq (2m)^{-1}$ for $i = 1, \dots, \tau$. Therefore
	\begin{displaymath}
	c^{\rm c}_\mathfrak{S}(k,p) \geq \frac{\big(|\{x \colon M^{\rm c}_k(f)(x) > 1/(3m)\}|\big)^{1/p}}{3m} \gtrsim \tau^{1/p} m^{1/p-1}. \eqno \qed 
	\end{displaymath} 
\end{prof*}

Next, we introduce the basic space $\mathfrak{T} = \mathfrak{T}_{\tau, d, m} = (Y, \rho, \mu)$ for fixed $\tau \in \mathbb{N}$, $1 < d \leq 3$ and $m > 1$. Let $Y = \{y_0, y^\circ_1, \dots, y^\circ_\tau,  y_1', \dots, y_\tau' \}$. We use some auxilliary symbols for certain subsets of $Y$: $Y^\circ = \{y^\circ_1, \dots, y^\circ_\tau \}$, $Y' = \{y_1', \dots, y_\tau'\}$ and $Y_i = \{y^\circ_i, y_i'\}$ for $i=1, \dots , \tau$. Define $\rho = \rho_d$ by the formula
\begin{displaymath}
\rho(x,y) = \left\{ \begin{array}{rl}
1 & \textrm{if } y_0 \in \{x,y\} \subset Y \setminus Y' \textrm{ or } \{x,y\} = Y_i, \ i \in \{1, \dots, \tau\}, \\
(d+1)/2 & \textrm{if } \{x,y\} \subset Y^\circ \textrm{ or } \{x,y\} \subset Y \setminus Y^\circ,\\
d & \textrm{otherwise,} \end{array} \right. 
\end{displaymath}
where $x$ and $y$ are two different elements of $Y$, and $\mu = \mu_m$ by letting $|\{y_0\}| = 1$, $|\{y^\circ_i\}| = 1 / \tau$ and $|\{y_i'\}| = m$ for $i=1, \dots, \tau$. Figure 3 shows a model of the space $\mathfrak{T}$. Adding an imaginary point at the top makes $\rho$ easily readable as a minor modification of the geodesic distance on the graph.
\begin{figure}[H]
	\begin{tikzpicture}
	[scale=.7,auto=left,every node/.style={circle,fill,inner sep = 2pt}]
	\node[label={[yshift=-1cm]$y_0$}] (n0) at (8,1) {};
	\node[label=$y^\circ_{1}$] (n1) at (6,3)  {};
	\node[label={[xshift=0.16cm]$y^\circ_{2}$}] (n2) at (7,3)  {};
	\node[label={[xshift=-0.24cm, yshift=-0.2cm]$y^\circ_{\tau-1}$}] (n3) at (9,3)  {};
	\node[label={[yshift=-0.05cm]$y^\circ_{\tau}$}] (n4) at (10,3)  {};
	\node[dots, scale=2] (n9) at (8,3)  {...};
	\node[label=$y_{1}'$] (n5) at (4,5)  {};
	\node[label={[xshift=-0.05cm]$y_2'$}] (n6) at (6,5)  {};
	\node[label={[xshift=0.3cm, yshift=-0.24cm]$y_{\tau-1}'$}] (n7) at (10,5)  {};
	\node[label={[yshift=-0.06cm]$y_{\tau}'$}] (n8) at (12,5)  {};
	\node[dots, scale=2] (n10) at (8,5)  {...};
	\node[label=$ $] (im) at (8,6.65)  {};
	\foreach \from/\to in {n0/n1, n0/n2, n0/n3, n0/n4, n1/n5, n2/n6, n3/n7, n4/n8}
	\draw (\from) -- (\to);
	
	
	\draw (6,5) arc (155:103.5:3);
	\draw (10,5) arc (25:76.5:3);
	\draw (4,5) arc (125:100:10);
	\draw (12,5) arc (55:80:10);
	\end{tikzpicture}
	\caption{The model of the space $\mathfrak{T}$.}
\end{figure}

Once again we explicitly describe any ball:
\begin{displaymath}
B(y_0,r) = \left\{ \begin{array}{rl}
\{y_0\} & \textrm{for } 0 < r \leq 1, \\
Y \setminus Y' & \textrm{for } 1 < r \leq (d+1)/2, \\
Y & \textrm{for } (d+1)/2 < r, \end{array} \right.
\end{displaymath} 
\noindent and, for $i \in \{1, \dots, \tau\}$,
\begin{displaymath}
B(y^\circ_i,r) = \left\{ \begin{array}{rl}
\{y^\circ_i\} & \textrm{for } 0 < r \leq 1, \\
\{y_0\} \cup Y_i & \textrm{for } 1 < r \leq (d+1)/2,  \\
\{y_0, y_i'\} \cup Y^\circ & \textrm{for } (d+1)/2 < r \leq d,  \\
Y & \textrm{for } d < r, \end{array} \right.
\end{displaymath}
and
\begin{displaymath}
B(y_i',r) = \left\{ \begin{array}{rl}
\{y_i'\} & \textrm{for } 0 < r \leq 1, \\
Y_i & \textrm{for } 1 < r \leq (d+1)/2,  \\
\{y_0, y^\circ_i\} \cup Y'  & \textrm{for } (d+1)/2 < r \leq d,  \\
Y & \textrm{for } d < r. \end{array} \right.
\end{displaymath} 

\begin{lemma}
	Let $\mathfrak{T}$ be the metric measure space defined as above. Then $c^{\rm c}_\mathfrak{T}(k,p) \simeq 1$ for $k \geq 1$ and $p \in [1, \infty]$, while
	\begin{displaymath}
	c_\mathfrak{T}(k,p) \simeq \left\{ \begin{array}{rl}
	\max\{1, \tau^{1/p}  m^{1/p-1}\} & \textrm{if } 1 \leq k < d \textrm{ and } p \in [1, \infty),  \\
	1 & \textrm{if } k \geq d \textrm{ or } p = \infty.\end{array} \right. 
	\end{displaymath}
\end{lemma}

\begin{prof*}
	First of all, it is easy to see that $c^{\rm c}_\mathfrak{T}(k,p) \geq 1$ and $c_\mathfrak{T}(k,p) \geq 1$ for any $k$ and $p$. Moreover, both $c^{\rm c}_\mathfrak{T}(k,p)$ and $c_\mathfrak{T}(k,p)$ are non-increasing as functions of $k$. Therefore, to prove that $c^{\rm c}_\mathfrak{T}(k,p) \simeq 1$, it suffices to show $c^{\rm c}_\mathfrak{T}(1,p) \lesssim 1$. Let $f \geq 0$ be a function on $Y$. Observe that $\max\{M^{\rm c}(f)(y)\colon y \in Y\} = \max\{f(y)\colon y \in Y\}$ which implies $c^{\rm c}_\mathfrak{T}(1,\infty) = 1$. Now assume that $p \in [1, \infty)$. We have
	\begin{displaymath}
	M^{\rm c}(f)(y_0) \leq \max\{f(y_0), A_{Y \setminus Y'}(f), A_Y(f)\}.
	\end{displaymath}
	Next, because $|\{y_0, y^\circ_i\} \cup Y'| \geq |Y'| \geq |Y|/3$, $i \in \{1, \dots, \tau\}$, then $A_{\{y_0, y^\circ_i\} \cup Y'}(f) \leq 3 A_Y(f)$ and hence
	\begin{displaymath}
	M^{\rm c}(f)(y_i') \leq \max\{f(y_i'), A_{Y_i}(f), 3 A_Y(f)\}.
	\end{displaymath}
	Finally, observe that
	\begin{displaymath}
	A_{\{y_0\} \cup Y_i}(f) \leq \max\{f(y_0), A_{Y_i}(f) \} \leq \max\{M^{\rm c}(f)(y_0), M^{\rm c}(f)(y_i') \},
	\end{displaymath}
	and
	\begin{displaymath}
	A_{\{y_0, y_i'\} \cup Y^\circ }(f) \leq \max\{A_{Y \setminus Y'}(f), f(y_i')\} \leq \max\{M^{\rm c}(f)(y_0), M^{\rm c}(f)(y_i') \},
	\end{displaymath}
	which implies
	\begin{displaymath}
	M^{\rm c}(f)(y^\circ_i) \leq \max\{f(y^\circ_i), M^{\rm c}(f)(y_0), M^{\rm c}(f)(y_i'), A_Y(f)\}.
	\end{displaymath}
	Since $|\{y^\circ_i\}| \leq |\{y_i'\}|$ and $\sum_{i=1}^{\tau} |\{y^\circ_i\}| = |\{y_0\}|$, we can estimate $\|M^{\rm c}(f)\|_p^p$ by
	\begin{displaymath}
	2 \, \Big( \sum_{y \in Y} f(y)^p \, |\{y\}| + 3^p \, A_Y(f)^p \, |Y| + A_{Y \setminus Y'}(f)^p \, |\{y_0\}| + \sum_{i=1}^{\tau} A_{Y_i}(f)^p \, |\{y_i'\}| \Big).
	\end{displaymath} 
	Applying H\"older's inequality we get $\|M^{\rm c}(f)\|_p^p \leq 2^{2p-1} \, (3^p + 3) \|f\|_p^p$ and thus $c^{\rm c}_\mathfrak{T}(1,p) \leq 24$.

	From now on we discuss only the non-centered case. It is easy to verify that $c_\mathfrak{T}(k,p) \simeq 1$ if $k \geq d$ or $ p = \infty$, by using exactly the same argument as in the proof of Lemma 4. In the next step we prove that $c_\mathfrak{T}(k,p) \lesssim \max\{1, \tau^{1/p}  m^{1/p-1}\}$ for $1 \leq k < d$ and $p \in [1, \infty)$. It suffices to show $c_\mathfrak{T}(1,p) \lesssim \max\{1, \tau^{1/p}  m^{1/p-1}\}$. Take $f \geq 0$ and observe that we have
	\begin{displaymath}
	M(f)(y_0) \leq \max\{f(y_0), A_{Y \setminus Y'}(f), M(f)(y_1'), \dots, M(f)(y_\tau')\},
	\end{displaymath}
	since $\{y_0\}$ and $Y \setminus Y'$ are the only balls containing $y_0$ and disjoint with $Y'$. Furthermore,
	\begin{displaymath}
	M(f)(y^\circ_i) \leq \max\{f(y^\circ_i), M(f)(y_0), M(f)(y_1'), \dots, M(f)(y_\tau')\},
	\end{displaymath}
	for any $i \in \{1, \dots, \tau\}$. Notice that if $y_i' \in B \subset Y$, then either $B \subset \{y_0, y_i'\} \cup Y^\circ$ or $|B| \geq |Y| / 3$. Since $|\{y_i'\}| \geq |\{y_0, y_i'\} \cup Y^\circ|/3$, we get
	\begin{displaymath}
	M(f)(y_i') \leq 3 \, \max\{A_{\{y_0, y_i'\} \cup Y^\circ}(f), A_Y(f)\}.
	\end{displaymath}
	Therefore we can estimate $\|M(f)\|_p^p$ by
	\begin{displaymath}
	3 \, \Big( \sum_{y \in Y \setminus Y'} f(y)^p \, |\{y\}| + A_{Y \setminus Y'}(f)^p \, |\{y_0\}| + 3^p  A_Y(f)^p \, |Y'| +  3^p  \sum_{i=1}^{\tau} A_{ \{y_0, y_i'\} \cup Y^\circ }(f)^p \, |\{y_i'\}| \Big).
	\end{displaymath}
	Since
	\begin{align*}
	\sum_{i=1}^{\tau} A_{ \{y_0, y_i'\} \cup Y^\circ }(f)^p \, |\{y_i'\}| & \leq 2^p \, \sum_{i=1}^{\tau} \Big( f(y_i')^p |\{y_i'\}|^p + \|f \cdot \chi_{Y \setminus Y'}\|_1^p \Big) |\{y_i'\}|^{1-p} \\
	& \leq 2^p \, \sum_{i=1}^\tau f(y_i')^p |\{y_i'\}| + 2^p \tau m^{1-p} \, \|f \cdot \chi_{Y \setminus Y'}\|_1^p,
	\end{align*}
	we can apply H\"older's inequality to get
	\begin{displaymath}
	\|M(f)\|_p^p \leq 3 \cdot 5^{p-1}  \Big( 1 + 1 + 3^p +6^p + 3^p \, 2^{2p-1} \tau m^{1-p} \Big) \|f\|_p^p,
	\end{displaymath}
	and, consequently, to obtain $c_\mathfrak{T}(1,p) \lesssim \max\{1, \tau^{1/p} m^{1/p-1}\}$.
	
	Finally, take $f = \delta_{\{y_0\}}$. If $1 \leq k < d$ and $p \in [1, \infty)$, then we have $\|f\|_p = 1$ and $M^{\rm c}_k(f)(y_i') = 1 / (1+ 1/\tau + m) \geq (3m)^{-1}$ for $i = 1, \dots, \tau$. Therefore
	\begin{displaymath}
	c^{\rm c}_\mathfrak{T}(k,p) \geq  \frac{\big(|\{ x \colon M^{\rm c}_k(f)(x) > 1/(4m)\}|\big)^{1/p}}{4m}  \gtrsim \tau^{1/p} m^{1/p-1}. \eqno \qed 
	\end{displaymath} 
\end{prof*}

\subsection{Auxiliary combinations of basic spaces}

Several times in Sections 3.2 and 3.3 we construct a metric measure space $\mathfrak{X}$ by using Proposition 1 for some $k_0$ and $\mathfrak{X}_n$, $n \in \Lambda$. It is worth noting here that, according to our notation, for $k \leq k_0$ and $p \in [1, \infty)$ we have
\begin{displaymath}
c_\mathfrak{X}(k,p) \gtrsim \sup\{c_{\mathfrak{X}_n}(k,p) \colon n \in \Lambda\}, \quad \Big(c_\mathfrak{X}(k,p)\Big)^p \lesssim \Big( \sup\{c_{\mathfrak{X}_n}(k,p) \colon n \in \Lambda\} \Big)^p + 1,
\end{displaymath}
which can be easily verified by following the proof of Proposition 1. Moreover, by taking $f \equiv 1$, one can easily see that $c_{\mathfrak{X}_n}(k,p) \geq 1$ for each $n \in \Lambda$. Combining these facts gives
\begin{displaymath}
c_\mathfrak{X}(k,p) \simeq \sup\{c_{\mathfrak{X}_n}(k,p) \colon n \in \Lambda\}, \quad k \leq k_0, \ p \in [1, \infty).
\end{displaymath}
An analogous statement for the centered operator also holds.  

\begin{lemma}
	Fix $1 \leq k < 2$, $p \in [1, \infty)$, $0 < \epsilon \leq 1/4$, $0 < \delta < 2-k$ and $ N \in \mathbb{N}$. For $n > N$ let $\mathfrak{X}_n = \mathfrak{S}_{\tau_n, d_n, m_n}$, where $\tau_n = N^{2p} {\lfloor}n^{p(p-1)/\epsilon}{\rfloor}$, $d_n = k + \delta / n$ and $m_n = n^{p/\epsilon}$. Define $\mathfrak{X}$ by using Proposition 1 for $k_0 = k + \delta$ and $\mathfrak{X}_n$, $n > N$. Then
	\begin{displaymath}
	c_\mathfrak{X}(k',p') \simeq c^{\rm c}_\mathfrak{X}(k',p'),
	\end{displaymath}
	for $k' \geq 1$ and $1 \leq p' \leq \infty$. Moreover
	\begin{enumerate}[label={\normalfont (\arabic*)}]
		\item $c_\mathfrak{X}(k',p') \simeq 1$ if $k' \geq k+\delta$ or $p' \geq p + 4\epsilon$,
		\item $c_\mathfrak{X}(k',p') < \infty$ if $k' > k$,
		\item $c_\mathfrak{X}(k',p') = \infty$ if $k' \leq k$ and $p' < p$,
		\item $c_\mathfrak{X}(k',p') \lesssim N^2$ if $k' \geq 1$ and $p' \geq p$,
		\item $c_\mathfrak{X}(k',p') \gtrsim N^{1/2} $ if $k' \leq k$ and $p' \in [p, p+\epsilon]$.  	
	\end{enumerate}	 
\end{lemma}

Here the symbol ${\lfloor} \cdot {\rfloor}$ refers to the floor function. Figure 4 describes the behavior of the function $c_\mathfrak{X}(k',p')$ (and thus also of $c^{\rm c}_\mathfrak{X}(k',p')$).

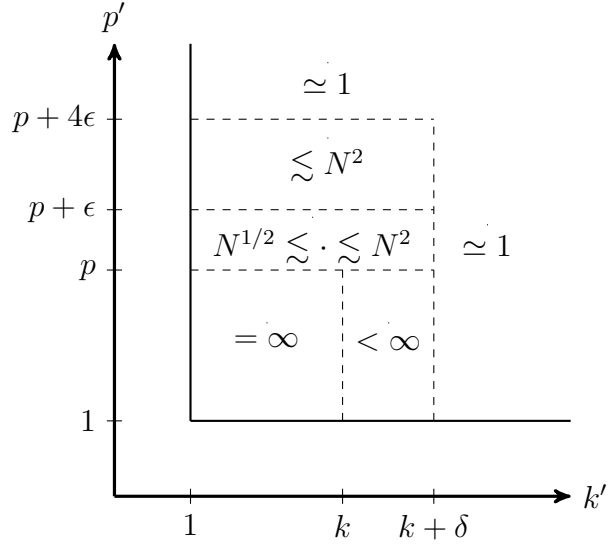
\begin{figure}[H]
	\begin{tikzpicture}[
	axis/.style={very thick, ->, >=stealth'},
	important line/.style={thick},
	dashed line/.style={dashed, thin},
	pile/.style={thick, ->, >=stealth', shorten <=2pt, shorten
		>=2pt},
	every node/.style={color=black}
	]
	\draw[axis] (0,0)  -- (6,0) node(xline)[right]
	{$k'$};
	\draw[axis] (0,0) -- (0,6) node(yline)[above] {$p'$};
	
	\draw[important line] (1,1) -- (6,1);
	\draw[important line] (1,1) -- (1,6);
	
	\draw[dashed line] (1,3) -- (4.2,3);
	\draw[dashed line] (1,3.8) -- (4.2,3.8);
	\draw[dashed line] (1,5) -- (4.2,5);
	
	\draw[dashed line] (3,1) -- (3,3);
	\draw[dashed line] (4.2,1) -- (4.2,5);		
	
	\draw (-0.1,1) node[left] {$1$} -- (0.1,1);
	\draw (-0.1,3) node[left] {$p$} -- (0.1,3);
	\draw (-0.1,3.8) node[left] {$p+\epsilon$} -- (0.1,3.8);
	\draw (-0.1,5) node[left] {$p+4\epsilon$} -- (0.1,5);
	
	\draw (1,-0.1) node[below] {$1$} -- (1,0.1);
	\draw (3,-0.1) node[below] {$k$} -- (3,0.1);
	\draw (4.2,-0.1) node[below] {$k+\delta$} -- (4.2,0.1);
	
	\draw (2,2.3) node[below] {$= \infty$} -- (2,2.3);
	\draw (3.6,2.3) node[below] {$< \infty$} -- (3.6,2.3);
	
	\draw (2.6,3.7) node[below] {$ N^{1/2} \lesssim \cdot \lesssim N^2 $} -- (2.6,3.7);
	
	\draw (4.9,3.6) node[below] {$ \simeq 1 $} -- (4.9,3.6);
	
	\draw (2.8,4.75) node[below] {$ \lesssim N^2 $} -- (2.8,4.75);
	
	\draw (2.8,5.75) node[below] {$ \simeq 1 $} -- (2.8,5.75);
	
	\end{tikzpicture}
	\caption{The behavior of the function $c_\mathfrak{X}(k',p')$.}
\end{figure}

\begin{prof*}
	Observe that $c_\mathfrak{X}(k',p') \simeq c^{\rm c}_\mathfrak{X}(k',p')$ for $k' \leq k_0$ and $1 \leq p' \leq \infty$, since 	
	\begin{displaymath}
	c_{\mathfrak{X}_n}(k',p') \simeq c^{\rm c}_{\mathfrak{X}_n}(k',p'), \qquad n > N,
	\end{displaymath}
	for that $k'$ and $p'$, and the same is true if we take the supremum over $n$. Moreover, $k_0 \geq d_n$ for each $n > N$, which implies $c_\mathfrak{X}(k_0, p') \simeq c^{\rm c}_\mathfrak{X}(k_0, p') \simeq 1$. Combining this with the fact that $c_\mathfrak{X}(k', p')$ and $c^{\rm c}_\mathfrak{X}(k', p')$ are estimated by $1$ from below and non-increasing as functions of $k'$ we conclude that $c_\mathfrak{X}(k',p') \simeq c^{\rm c}_\mathfrak{X}(k',p')$ holds for the full range of the parameters $k'$ and $p'$. Now, to prove $(1)$ it suffices to show that $c_\mathfrak{X}(1,p') \simeq 1$ for $p + 4\epsilon \leq p' < \infty$. For that $p'$ and $n > N$ we have the following inequality
	\begin{align*}
	c_{\mathfrak{X}_n}(1, p') & \lesssim 1 + N^{2p/p'} \cdot n^{p(p-1)/(\epsilon p')} \cdot n^{p(1-p')/(\epsilon p')} \\
	& \lesssim 1 + N^{2p/p'} \cdot n^{p(p-p')/(\epsilon p')} \lesssim 1 + N^{2} n^{-2} \lesssim 1,
	\end{align*}
	since $1 \leq p/p' \leq 2$ and $p-p' \leq - 4 \epsilon$. This implies 
	\begin{displaymath}
	c_\mathfrak{X}(1,p') \lesssim \sup\{c_{\mathfrak{X}_n}(1, p') \colon n > N\} \lesssim 1.
	\end{displaymath}
	Condition $(2)$, in turn, is a simple consequence of the fact that $d_n > k'$ only for finitely many $n$ if $k' > k$. Next, take $k' \geq k$ and $p' < p$ (we can do this only if $p \neq 1$). Then
	\begin{align*}
	c_\mathfrak{X}(k', p') \gtrsim \limsup_{n \rightarrow \infty} c_{\mathfrak{X}_n}(k', p') & \gtrsim \lim_{n \rightarrow \infty} N^{2p/p'}  \cdot n^{p(p-1)/(\epsilon p')} \cdot n^{p(1-p')/(\epsilon p')} \\
	& \gtrsim \lim_{n \rightarrow \infty} N^{2p/p'} \cdot n^{p(p-p')/(\epsilon p')} = \infty,
	\end{align*}
	and $(3)$ holds. To prove $(4)$ assume that $p' \geq p$. For each $n > N$ we have
	\begin{align*}
	c_{\mathfrak{X}_n}(k', p') & \lesssim 1 + N^{2p/p'} \cdot n^{p(p-1)/(\epsilon p')} \cdot n^{p(1-p')/(\epsilon p')} \\
	& \lesssim 1+ N^{2p/p'} \cdot n^{p(p-p')/(\epsilon p')} \leq 1 + N^{2} \cdot 1 \lesssim N^{2},
	\end{align*}
	and therefore
	\begin{displaymath}
	c_\mathfrak{X}(k', p') \lesssim \sup\{ c_{\mathfrak{X}_n}(k', p') \colon n > N\} \lesssim N^2.
	\end{displaymath}
	Finally, take $k' \geq k$ and $p' \in [p, p+\epsilon]$. Since $3/4 \leq p/p' \leq 1$ and $-\epsilon \leq p-p' \leq 0$, we have
	\begin{align*}
	c_\mathfrak{X}(k', p') \gtrsim c_{\mathfrak{X}_{2N}}(k', p') & \gtrsim N^{2p/p'} \cdot (2N)^{p(p-1)/(\epsilon p')} \cdot (2N)^{p(1-p')/(\epsilon p')} \\
	& \gtrsim  N^{2p/p'} \cdot (2N)^{p(p-p')/(\epsilon p')} \gtrsim N^{3/2} \cdot N^{-1} = N^{1/2},
	\end{align*}
	which justifies $(5)$ and completes the proof. $\raggedright \hfill \qed$	
\end{prof*}

\begin{lemma}
	Fix $1 < k \leq 2$ (respectively, $1 \leq k < 2$) and let $\mathfrak{X}_n = \mathfrak{S}_{\tau_n, d_n, m_n}$ with $\tau_n = n$, $d_n = k$ (respectively, $d_n = k + \frac{2-k}{n}$) and $m_n = 2$. Define $\mathfrak{X}$ by using Proposition 1 for $k_0 = 2$ and $\mathfrak{X}_n$, $n \in \mathbb{N}$. Then $c_{\mathfrak{X}}(k',p) = \infty$ if and only if $k'<k$ (respectively, $k' \leq k$) and $p \neq \infty$, and the same is true for $c^{\rm c}_{\mathfrak{X}}(k',p)$ replacing $c_{\mathfrak{X}}(k',p)$.
\end{lemma}
 \begin{prof*}
 	We prove only the first case (the second one can be obtained very similarly). Assume that $p \neq \infty$ since the case $p = \infty$ is trivial. If $k' < k$, then for any $n \in \mathbb{N}$ we have $k' < d_n$ and hence $c^{\rm c}_{\mathfrak{X}_n}(k',p) \simeq n^{1/p}$. Therefore 
 	\begin{displaymath}
 	c_{\mathfrak{X}}(k',p) \geq c^{\rm c}_{\mathfrak{X}}(k',p) \gtrsim \limsup_{n \rightarrow \infty} c_{\mathfrak{X}_n}(k',p) \simeq \lim_{n \rightarrow \infty} n^{1/p} = \infty.
 	\end{displaymath}
 	In turn, if $k' \geq k = d_n$, then $c_{\mathfrak{X}_n}(k,p) \simeq 1$ gives 
 	\begin{displaymath}
 	c^{\rm c}_{\mathfrak{X}}(k',p) \leq c_{\mathfrak{X}}(k',p) \leq c_{\mathfrak{X}}(k,p) \lesssim 1 < \infty. \eqno \qed
 	\end{displaymath}
 \end{prof*}
 
 We notice (without furnishing the detailed proof) that one can obtain the following counterparts of Lemmas 6 and 7 by using the adequate spaces $\mathfrak{T}_{\tau_n, d_n, m_n}$ instead of $\mathfrak{S}_{\tau_n, d_n, m_n}$.
 
 \begin{lemma6}
 	Fix $1 \leq k < 3$, $p \in [1, \infty)$, $0 < \epsilon \leq 1/4$, $0 < \delta < 3-k$ and $ N \in \mathbb{N}$. For $n > N$ let $\mathfrak{Y}_n = \mathfrak{T}_{\tau_n, d_n, m_n}$, where $\tau_n = N^{2p} {\lfloor}n^{p(p-1)/\epsilon}{\rfloor}$, $d_n = k + \delta / n$ and $m_n = n^{p/\epsilon}$. Define $\mathfrak{Y}$ by using Proposition 1 for $k_0 = k + \delta$ and $\mathfrak{Y}_n$, $n > N$. Then
 	\begin{itemize}
 		\item $c^{\rm c}_\mathfrak{Y}(k',p') \simeq 1$ for $k' \geq 1$ and $1 \leq p \leq \infty$,
 		\item conditions $(1) - (5)$ from Lemma 6 hold with $c_\mathfrak{Y}(k',p')$ replacing $c_\mathfrak{X}(k',p')$.  	
 	\end{itemize}	 
 \end{lemma6}
 
 \begin{lemma7}
 	Fix $1 < k \leq 3$ (respectively, $1 \leq k < 3$) and let $\mathfrak{Y}_n = \mathfrak{T}_{\tau_n, d_n, m_n}$ with $\tau_n = n$, $d_n = k$ (respectively, $d_n = k + \frac{3-k}{n}$) and $m_n = 2$. Define $\mathfrak{Y}$ by using Proposition 5 for $k_0 = 3$ and $\mathfrak{Y}_n$, $n \in \mathbb{N}$. Then $c_{\mathfrak{Y}}(k',p) = \infty$ if and only if $k'<k$ (respectively, $k' \leq k$) and $p \neq \infty$, and the same is true for $c^{\rm c}_{\mathfrak{Y}}(k',p)$ replacing $c_{\mathfrak{Y}}(k',p)$.
 \end{lemma7}
 
 As the last thing in this section let us point out here that each of the spaces constructed by using Lemmas $6$, $7$, $6'$ or $7'$ is non-doubling and hence the same will be true for the spaces constructed in the proof of Theorem $2$. 
 
 \subsection{Proof of Theorem 2}
 
 	In the first step we construct a metric measure space $\mathfrak{X}$ for which the associated modified maximal operators $M_k^{\rm c}$ and $M_k$, $k \in [1, 2)$, are of weak type $(p,p)$ if and only if $p > h^{\rm c}(k)$ or $p = \infty$, while $M_2$ is of weak type $(1,1)$. The last property can be easily checked, since the basic spaces $\mathfrak{S}$ will be used to build $\mathfrak{X}$. Consider the case $h^{\rm c}(k) < \infty$ for each $k \in [1,2]$. Let us introduce a countable set $\Sigma = \Sigma_1 \cup \Sigma_2 = \{k_1, k_2, \dots\}$, where $\Sigma_1$ is the set of all $k \in [1,2)$ for which $\lim_{k' \rightarrow k^+} h^{\rm c}(k') < h^{\rm c}(k)$ (the case $\Sigma_1 = \emptyset$ is possible) and $\Sigma_2$ is a dense subset of $(1,2)$ that has no common points with $\Sigma_1$. By induction we will construct a family of metric measure spaces $\mathfrak{X}_{n,m}$, $n, m \in \mathbb{N}$, and then we will obtain $\mathfrak{X}$ by using Proposition 1.
 	
 	Take $k_1 \in \Sigma$ and let $0 < \delta_1 < 2-k_1$ be such that
 	\begin{displaymath}
 	h^{\rm c}(k') \geq \lim_{k \rightarrow k_1^+} h^{\rm c}(k) - 1, \quad k' \leq k_1 + \delta_1.
 	\end{displaymath}
 	For each $m \in \mathbb{N}$ denote by $\mathfrak{X}_{1,m}$ the space constructed by using Lemma 6 with $k=k_1$, $p=h^{\rm c}(k_1)$, $\epsilon = 1/(4m)$, $\delta = \delta_1 / m$ and $N=m$. Now, let $n \geq 2$ and suppose that for each $j \in \{1, \dots, n-1\}$ and $m \in \mathbb{N}$ the space $\mathfrak{X}_{j,m}$ has been already constructed. We choose $0 < \delta_n < 2-k_n$ such that the following conditions are satisfied
 	\begin{itemize}
 		\item $h^{\rm c}(k') \geq \lim_{k \rightarrow k_n^+} h^{\rm c}(k) - 1/n$ for $k' \leq k_n + \delta_n$,
 		\item if $k_j > k_n$ for some $j \in \{1, \dots, n-1\}$, then $k_n + \delta_n < k_j$.
 	\end{itemize} 
 	For each $m \in \mathbb{N}$ we construct $\mathfrak{X}_{n,m}$ as in Lemma 6 with $k=k_n$, $p=h^{\rm c}(k_n)$, $\epsilon = 1/(4m)$, $\delta = \delta_n / m$ and $N=m$. Finally, denote by $\mathfrak{X}$ the space constructed by using Proposition 1 with $k_0 = 2$ for $\mathfrak{X}_{n,m}$, $n, m \in \mathbb{N}$. It suffices to show that for each $k \in [1,2)$ we have $c^{\rm c}_{\mathfrak{X}}(k, h^{\rm c}(k)) = \infty$, while $c_{\mathfrak{X}}(k, p) < \infty$ if $p > h^{\rm c}(k)$. 
 	
 	Fix $k \in [1,2)$. If $\lim_{k' \rightarrow k^+} h^{\rm c}(k') < h^{\rm c}(k)$, then $k = k_{n} \in \Sigma$ for some $n \in \mathbb{N}$. Therefore $c^{\rm c}_{\mathfrak{X}_{n,m}}(k, h^{\rm c}(k)) \gtrsim m^{1/2}$ which implies
 	\begin{displaymath}
 	c^{\rm c}_{\mathfrak{X}}(k, h^{\rm c}(k)) \gtrsim \sup\{ c^{\rm c}_{\mathfrak{X}_{n,m}}(k, h^{\rm c}(k)) \colon m \in \mathbb{N} \} \gtrsim \lim_{m \rightarrow \infty} m^{1/2} = \infty.
 	\end{displaymath}
 	In turn, if $\lim_{k' \rightarrow k^+} h^{\rm c}(k') = h^{\rm c}(k)$, then for any $j = 1, 2, \dots$ we can choose a point $k_{n_j} \in \Sigma$ such that $k_{n_j} > k$ and $h^{\rm c}(k_{n_j}) > h^{\rm c}(k) - 1/(4j)$. Hence $c^{\rm c}_{\mathfrak{X}}(k, h^{\rm c}(k)) \gtrsim c^{\rm c}_{\mathfrak{X}_{n_j,j}}(k, h^{\rm c}(k)) \gtrsim j^{1/2}$ and letting $j \rightarrow \infty$ we obtain $c^{\rm c}_{\mathfrak{X}}(k, h^{\rm c}(k)) = \infty$. 
 	
 	Next, fix $k \in [1,2)$ and let $h^{\rm c}(k) < p < h^{\rm c}(k) + 1$. It is obvious that $c_{\mathfrak{X}_{n,m}}(k,p) < \infty$ for any fixed $n$ and $m$. We will prove that
 	\begin{displaymath}
 	\sup\{c_{\mathfrak{X}_{n,m}}(k,p) \colon n,m \in \mathbb{N} \} < \infty.
 	\end{displaymath}
 	Let $n_0$ be such that
 	\begin{displaymath}
 	h^{\rm c}(k) + \frac{1}{n_0+1} \leq p < h^{\rm c}(k) + \frac{1}{n_0}.
 	\end{displaymath} 
 	Take $n \in \mathbb{N}$ such that $k \notin [k_n, k_n + \delta_n)$. With this assumption we obtain $c_{\mathfrak{X}_{n,m}}(k,p) \simeq 1$ for $m \geq n_0+1$. In turn, if $m < n_0+1$, then $c_{\mathfrak{X}_{n,m}}(k,p) \lesssim m^2 \leq (n_0 + 1)^2$. Next, let $n \in \mathbb{N}$ be such that $k \in [k_n, k_n + \delta_n)$. There exists $m_0 = m_0(n)$ such that $k \notin [k_n, k_n + \delta_{n,m_0})$ or $h^{\rm c}(k_n) + 1/m_0 < p$. This implies $c_{\mathfrak{X}_{n,m}}(k,p) \simeq 1$ for any $m \geq m_0$. Hence, we deduce that $\sup\{c_{\mathfrak{X}_{n,m}}(k,p) \colon n,m \in \mathbb{N} \} < \infty$, if there is a finite number of $n$ such that $k \in [k_n, k_n + \delta_{n})$. Otherwise, choose $l \geq 2(n_0+1)$ such that $k \in [k_l, k_l + \delta_{l})$. If $k \in [k_n, k_n + \delta_{n})$ for some $n > l$, then
 	\begin{displaymath}
 	h^{\rm c}(k) \geq \lim_{k' \rightarrow k_l^+} h^{\rm c}(k') - \frac{1}{l} \geq h^{\rm c}(k_n) - \frac{1}{2(n_0+1)},
 	\end{displaymath} 
 	since $k_n \in (k_l, k_l + \delta_{l})$, which implies
 	\begin{displaymath}
 	p \geq h^{\rm c}(k) + \frac{1}{n_0+1} \geq h^{\rm c}(k_n) + \frac{1}{2(n_0+1)}.
 	\end{displaymath}
 	Hence, for that $n$, if $m \geq 2(n_0+1)$, then $c_{\mathfrak{X}_{n,m}}(k,p) \simeq 1$. Since $c_{\mathfrak{X}_{n,m}}(k,p) \lesssim 4(n_0+1)^2$ for $m < 2(n_0+1)$, we conclude that $\sup\{c_{\mathfrak{X}_{n,m}}(k,p) \colon n,m \in \mathbb{N} \} < \infty.$
 	
 	Suppose now that $h^{\rm c}$ takes the value $\infty$ and denote $a = \sup\{k \colon h^{\rm c}(k) = \infty\}$. If $a = 2$, then we use the appropriate version of Lemma 7 with $k=2$ to choose $\mathfrak{X}$. Assume that $a < 2$. If $\lim_{k \rightarrow a^+} h^{\rm c}(k) = \infty$, then $h^{\rm c}$ is continuous at $a$ and we just construct $\mathfrak{X}$ in the same way as we did in the case $h^{\rm c} < \infty$, but now using $[a,2)$ and $(a,2)$ instead of $[1,2)$ and $(1,2)$, respectively. It is not hard to verify that $\mathfrak{X}$ has the expected properties. Otherwise, we introduce an auxiliary function $h'$ defined by the formula
 	\begin{displaymath}
 	h'(k) = \left\{ \begin{array}{rl}
 	h_0 & \textrm{if } 1 \leq k \leq a,   \\
 	h^{\rm c}(k) & \textrm{if } a < k \leq 2, \end{array} \right. 
 	\end{displaymath} 
 	where $h_0 = h^{\rm c}(a)$ if $h^{\rm c}(a) < \infty$ or $h_0 = \lim_{k \rightarrow a^+} h^{\rm c}(k)$ in the opposite case. Let $\mathfrak{X}'$ be the space constructed as before with $h'$ instead of $h^{\rm c}$. We use Proposition 1 with $k_0 = 2$ one more time and obtain $\mathfrak{X}$ combining $\mathfrak{X}'$ with the space from Lemma 7 with $k = a$ (we use the appropriate version of Lemma 7 depending on whether $h^{\rm c}(a) < \infty$ or $h^{\rm c}(a) = \infty$).
 	
 	In the second step we construct a metric measure space $\mathfrak{Y}$ for which the associated modified maximal operators $M_k$, $k \in [1, 3)$, are of weak type $(p,p)$ if and only if $p > h(k)$ or $p = \infty$, while the operators $M_k^{\rm c}$, $k \in [1,2)$, are of weak type $(p,p)$ for any $p \geq 1$. The method is similar to that which was used to construct $\mathfrak{X}$. The key point is to use Lemmas $6'$ and $7'$ instead of Lemmas 6 and 7, respectively, and to use Proposition 1 with $k_0 = 3$. We skip the technical details here.
 	
 	Finally, we build the metric measure space $\mathfrak{Z}$ by using Proposition 1 with $k_0 = 3$ for $\mathfrak{X}$ and $\mathfrak{Y}$. It is not hard to see that we have
 	\begin{displaymath}
 	\big( c^{\rm c}_{\mathfrak{Z}}(k,p) < \infty \big) \Leftrightarrow \big( \max\{c^{\rm c}_{\mathfrak{X}}(k,p), c^{\rm c}_{\mathfrak{Y}}(k,p)\} < \infty \big) \Leftrightarrow \big( p \in (h^{\rm c}(k), \infty] \big) , \quad k \in [1,2),
 	\end{displaymath}
 	and
 	\begin{displaymath}
 	\big( c_{\mathfrak{Z}}(k,p) < \infty \big) \Leftrightarrow \big( \max\{c_{\mathfrak{X}}(k,p), c_{\mathfrak{Y}}(k,p)\} < \infty \big) \Leftrightarrow \big( p \in (h(k), \infty] \big) , \quad k \in [1,3),
 	\end{displaymath}
 	(here we use the convention $(\infty, \infty] = \{\infty\}$) and hence $\mathfrak{Z}$ may be chosen as the metric measure space satisfying all the expected conditions. Thus, the proof of Theorem 2 is complete.

\subsection{Additional examples}

In the last section, for simplicity, we deal only with the centered operators. We will try to give a general idea of how the situation changes when we want $P_k^{\rm c}$, $k \in [1,2]$, to be of the form $[h^{\rm c}(k), \infty]$, where $h^{\rm c}$ is a fixed function.   

\begin{example}
	Let $h^{\rm c} \colon [1,2] \rightarrow [1, \infty]$ be a continuous non-increasing function with $h^{\rm c}(2)=1$. Then there exists a metric measure space $\mathfrak{Z}$ such that the associated modified maximal operators $M_k^{\rm c}$, $k \in [1, 2)$, are of weak type $(p,p)$ if and only if $p \geq h^{\rm c}(k)$. 
\end{example}

\begin{proof}
	If $h^{\rm c}(1) = 1$, then the result is trivial since $\mathfrak{Z}$ may be chosen to be $\{a\}$, the set of one point, equipped with the unique metric and counting measure. From now on assume that $h^{\rm c}(1) > 1$. Let us introduce an auxiliary set
\begin{displaymath}
\Omega = \{(k,p) \in \big( [1,2] \cap \mathbb{Q} \big) \times \big( [1, \infty) \cap \mathbb{Q} \big) \colon p < h^{\rm c}(k) \} = \{(k_n, p_n) \colon n \in \mathbb{N}\}.
\end{displaymath}	
For each $n \in \mathbb{N}$ we choose $0 < \delta_n < 2-k_n$ such that $p_n < h^{\rm c}(k_n + \delta_n)$. Denote by $\mathfrak{X}_n$ the space constructed by using Lemma 6 with $k = k_n$, $p = p_n$, $\epsilon = 1/4$, $\delta = \delta_n$ and $N = 1$. Then it is easy to show that $\mathfrak{Z}$ may be chosen to be the space constructed by using Proposition 1 with $k_0 = 2$ for $\mathfrak{X}_n$, $n \in \mathbb{N}$.
\end{proof}

The second example is more general in the sense that we take into account all metric measure spaces, not just those that satisfy the assumptions made at the beginning of the article. In the proof we will apply the estimates of the operator norm obtained by interpolation (see \cite{DB}, Theorem VIII.9.1, p. 392, for example). Moreover, we will use the basic fact that for any metric measure space $\mathfrak{X}$ we have
\begin{displaymath}
\lim_{k \rightarrow k_0^+} c^{\rm c}_{\mathfrak{X}}(k, p_0) = c^{\rm c}_{\mathfrak{X}}(k_0, p_0), \qquad k_0 \geq 1, \ p_0 \geq 1.
\end{displaymath} 

\begin{example}
Let $\mathbb{Q} \cap (1,2) = \{q_1, q_2, \dots\}$ and define
\begin{displaymath}
h^{\rm c}(k) = 2 - \sum_{i \colon q_i < k} \frac{1}{2^i}, \qquad k \in [1,2].
\end{displaymath}
Then there is no metric measure space such that for each $k \in [1,2]$ the associated maximal operator $M_k^{\rm c}$ is of weak type $(p,p)$ if and only if $p \geq h^{\rm c}(k)$. 
\end{example}

\begin{proof}
Suppose by contradiction that $\mathfrak{X}$ is such a space. First we show that for any $1 \leq a < b \leq 2$ and $N \in \mathbb{N}$ we can find $ a \leq a' < b' \leq b$ such that
\begin{displaymath}
c_{\mathfrak{X}}^{\rm c}(k, h^{\rm c}(k)) \geq N, \qquad k \in [a', b'].
\end{displaymath}
Indeed, take $q_i \in (a,b)$ and observe that
\begin{equation}\label{1}
\lim_{k \rightarrow q_i^+} c^{\rm c}_{\mathfrak{X}}\big(k, h^{\rm c}(q_i) - \frac{1}{2^{i+1}}\big) = \infty.
\end{equation}
Next, let $0 < \epsilon < 2 - q_i$. From the definition we have $h^{\rm c}(q_i)-\frac{1}{2^{i+1}} - h^{\rm c}(q_i + \epsilon) \geq \frac{1}{2^{i+1}}$. Moreover, note that $q_i \notin (1, q_i)$ implies $1 \leq h^{\rm c}(q_i) - \frac{1}{2^{i+1}} \leq 2$. Thus, if $c^{\rm c}_{\mathfrak{X}}(q_i+\epsilon, h^{\rm c}(q_i+\epsilon)) \leq N$, then by interpolation we obtain
\begin{displaymath}
c^{\rm c}_{\mathfrak{X}}(q_i+\epsilon, h^{\rm c}(q_i)-\frac{1}{2^{i+1}}) \leq 2 \Big( \frac{h^{\rm c}(q_i)-\frac{1}{2^{i+1}}}{h^{\rm c}(q_i)-\frac{1}{2^{i+1}} - h^{\rm c}(q_i + \epsilon)} \Big)^{\frac{1}{h^{\rm c}(q_i) - 1 /2^{i+1}}} N ^{1 - \frac{h^{\rm c}(q_i + \epsilon)}{h^{\rm c}(q_i) - 1 / 2^{i+1}}} \leq 2^{i+3}N.
\end{displaymath} 
Of course, in view of \eqref{1}, such an estimate cannot occur for sufficiently small values of $\epsilon$. Therefore, we can choose an interval $[a', b'] \subset (q_i, b] \subset [a,b]$ with the expected property. The rest of the proof consists of constructing inductively the sequence of closed intervals $[1,2] \supset [a_1, b_1] \supset [a_2, b_2] \supset \dots$, such that for each $n \in \mathbb{N}$ and $k \in [a_n, b_n]$ we have $c_{\mathfrak{X}}^{\rm c}(k, h^{\rm c}(k)) \geq n$. Clearly, $\bigcap_{n=1}^\infty [a_n, b_n] \neq \emptyset$ and $c_{\mathfrak{X}}^{\rm c}(k, h^{\rm c}(k)) = \infty$ for any $k \in \bigcap_{n=1}^\infty [a_n, b_n]$, which contradicts the assumption that $M_k^{\rm c}$ is of weak type $(h^{\rm c}(k),h^{\rm c}(k))$.
\end{proof}

\section*{Acknowledgement}
I would like to express my deep gratitude to my supervisor Professor Krzysztof Stempak for his suggestion to study the problem discussed in this article. I thank him for all the remarks made during the preparation of the manuscript. 

I am also indebted to the referee for a very thorough reading of the article and for constructive comments and hints resulting in an improvement of the presentation. 

Research was supported by the National Science Centre of Poland, project no. \linebreak 2016/21/N/ST1/01496.

\end{document}